\RequirePackage{fix-cm}
\documentclass[natbib]{svjour3}       
\smartqed  
\usepackage{amsmath}
\usepackage{mathrsfs}
\usepackage{amsfonts}
\usepackage{booktabs}
\usepackage{caption}
\usepackage{multirow}
\usepackage{rotating}
\usepackage{longtable}
 \usepackage{algorithm}
\usepackage{algorithmic}
\usepackage{graphicx,psfrag,color,pstricks,pst-grad}
\begin{document}
 
\title{Determination of the doubly-symmetric 
	periodic orbits in the 
	restricted three-body problem 
	and Hill's lunar problem
\thanks{Xu is supported by the National Nature Science Foundation of China (NSFC, Grant No. 11703006).}
}

\titlerunning{Determination of doubly-symmetric periodic orbits}        

\author{Xingbo Xu          
}


\institute{ Xingbo Xu \at
              Huaiyin Institute of Technology, Faculty of Mathematics and Physics,
              department of math.    \\
              Tel.: +86-13705237961\\
              \email{xuxingbo25@hotmail.com}  
}

\date{Received: 23 August 2022 / Accepted: 26 August 2022}

\maketitle

\begin{abstract}
  We review some recent progress on the research
  of the periodic orbits of the N-body problem,
  and propose a numerical scheme to determine
  the spatial doubly-symmetric periodic orbits (\textsc{SDSP}s for short).
  Both comet- and lunar-type \textsc{SDSP}s 
  in the circular restricted three-body problem 
  are computed, as well as the Hill-type \textsc{SDSP}s 
  in Hill's lunar problem.  
  Doubly symmetries are exploited so that the 
  \textsc{SDSP}s can be computed efficiently.
  The monodromy matrix can be calculated by the information
  of one fourth period.
  The periodicity conditions are solved by
  Broyden's method with a line-search, and the 
  algorithm is reviewed. Some numerical examples
  show that the scheme is very efficient. 
   For a fixed period ratio and a given acute angle,
   there exist sixteen cases of initial values.  For the restricted 
   three-body problem, the cases of ``Copenhagen problem'' and the 
   Sun-Jupiter-asteroid model are considered. 
   New \textsc{SDSP}s are also numerically found in Hill's lunar problem. 
   Though the period ratio should be small theoretically,
   some new periodic orbits are found when
   the ratio is not too small, 
   and most of the searched \textsc{SDSP}s
   are linearly stable.  
\end{abstract}

\section{Introduction}
\label{intro}
  
  
  Many versions of the n-body problems are studied within the theme
   of periodic orbits or KAM tori, as
   periodic orbits and quasi-periodic orbits are important in
   celestial mechanics and very useful 
   in understanding the motions of celestial bodies. 
   For the aspects of the applications of the KAM theories on celestial
   mechanics, one may refer to \citet{Celletti06,Biasco08,Meyer2011,
   ZhaoL14} and the references therein.
   Here, we mainly discuss the periodic orbits.
   For the non-integrable dynamical systems, 
   periodic orbits determine the skeleton of the phase space, and 
   stable ones usually accumulate onto the KAM tori. 
   
   The problem on the periodic orbits of the three-body problem
    is intractable ever since \citet{Poincare}. 
     Symmetric periodic orbits are prior to be studied, as 
   symmetries can usually reduce the complexity.
   Some earlier work on the periodic orbits can be referred to
    \citet{Hadji84}.  Three classical families of    
    periodic orbits of the three-body problem are of Euler-Lagrange type \citep{Meyer86, Sicardy10, HuSun10, Zhou17}, 
    figure-eight type \citep{ChencinerM, MunozA07, ChenKC09, HuSun09, GalanV14, YuG17} and ``Broucke–H\'{e}non'',``Schubart-like'' type 
    \citep{Henon74, Henon76, Martinez13, Ortega16, Voyatzis18, Kuang19}.      
    It is summarized in \citet{LiXM19} that only these three families and
    a few more new periodic orbits are numerically calculated before 
    a series of their work, and more than $2000$
    families of periodic orbits are found by taking use of the ``clean
    numerical simulation” (CNS) technique, the grid search
    and the Newton-Raphson method on a supercomputer.
    The CNS is based on the arbitrary order of Taylor series method
    with a convergence check \citep{Jorba05, Liao19}. 
    The periodicity conditions are based on the difference of two solutions 
    at different times, and the permissible error 
    is set to be less than $10^{-6}$.

   The \emph{restricted three-body problem} (\textsc{RTBP}) and Hill's
   lunar problem are well known as two basic mathematical models,
   in which many families of periodic orbits 
   have been shown to exist, analysed and calculated. 
    Periodic orbits of a nearly integrable Hamiltonian system 
    with a small parameter are usually shown to exist by 
    Poincar\'{e}'s continuation method or Arenstorf's implicit 
    function theorem, combinating several techniques such as
    canonical transforms, symplectic scaling,
    averaging, symmetry reduction and so on \citep{MeyerHO,CorsJM05}.
    If there exist orbits very far away from the primaries,
  these orbits are classified as the \emph{comet-type}. 
   If there exist orbits around and very close to one primary,
   these orbits are of \emph{Hill-type}. 
    \citet{Howison00I,Howison00II} showed the existence of 
     both the comet- and Hill-type
     \emph{spatial doubly-symmetric periodic orbits}
      (\textsc{SDSP}s) in the circular \textsc{RTBP} (\textsc{CRTBP}) 
      and the Hill-type \textsc{SDSP}s in Hill's lunar problem.    
      Howison and Meyer's results were generalized 
    to the cases of the restricted (n+1)-body problem with general
    homogeneous potential successively by \citet{LlibreR09}
    and \citet{LlibreS11},  where $n$ bodies form a regular polygon 
    or a nested central configuration.       
    \citet{XbXu19,XuHill20} showed that there exists a class of 
    Hill-type \textsc{SDSP}s around one oblate primary 
    in a generalized CRTBP and in Hill's lunar problem. 
    According to the known results \citep{Palacian06,Meyer2011},
    there exist spatial comet-type KAM tori in the elliptic
    \textsc{RTBP}.
    For the lunar three-body problem, which is also called the
    hierarchal triple system, \citet{ZhaoL14} showed that there
    exist families of periodic orbits accumulating onto several
    KAM tori by applying F\'{e}joz's finer version of KAM theorem.
    These analytical results reveal that some \textsc{SDSP}s exist and
    stable when the small perturbation parameter is sufficiently small.

     In \citet{BroerZhao17}, only one linearly stable 
   family is shown to exist though there are $16$ families of 
   Keplerian periodic orbits of De Sitter. This enlightens us to  
  numerically study Howison and Meyer's \textsc{SDSP}s 
   with different initial values. 
   If the initial values of \textsc{SDSP}s of the full system 
   can be numerically continued from those of the circular \textsc{SDSP}s 
   of the approximate system, 
   the \textsc{SDSP}s of the full system 
   will be of multiple revolutions and nearly circular.
   According to Poincar\'{e}'s classification, 
   these orbits in the \textsc{CRTBP} can be classified
   as the third sort of the first species, while these orbits in Hill's
   lunar problem can be classified as the third sort of the third species.
    However, we don't know the cases when the perturbation is
    not too small, thus the numerical work on
    these periodic orbits is necessary and interesting.
    The purpose of this paper is to study
    the numerical calulation and stability of the \textsc{SDSP}s in the \textsc{CRTBP} and Hill's lunar problem.

    The readers who are interested in the related history 
    and more aspects are encouraged to refer to 
    \citet{Henon69, Henon97, Henon03, BrunoV06, UrsOtto} 
    and the citations. 
    In the earlier stage, solutions of the periodicity conditions 
    are continued based on Newton's method. \citet{Macris75}
    calculated four families of \textsc{SDSP}s of the elliptic RTBP
    with one fixed mass $\mu=0.4$ and
    the eccentricity as the small parameter, 
    but all of the $107$ orbits are unstable.
    With the mass ratio of the Sun-Jupiter case,  
    \citet{Kazantis79} calculated $7$ families of \textsc{SDSP}s which
    bifurcate from the vertical-critical symmetric families
  \emph{a}, \emph{g}$_{1}$, \emph{g}$_{2}$, \emph{h}, \emph{i}, \emph{l} and \emph{m}. 
  Family \emph{a} goes from the Euler collinear libration point $L_{2}$. 
  Families \emph{g}$_{1}$ and \emph{g}$_{2}$
  are branches of the family \emph{g} which begins as
  direct circular orbits of Hill-type around the Jupiter.
  Families \emph{h} and \emph{i} begin from the retrograde and direct circular orbits of Hill-type around the Sun, respectively.       
  Families \emph{l} and \emph{m} are retrograde
  nearly circular periodic orbits of comet-type with positive and negative Jacobi constants, respectively. 
  In addition, families \emph{b} and \emph{c}
  begin from the infinitesimal orbits around 
  $L_{3}$ and $L_{1}$ respectively. Family \emph{f}
  begins from the infinitesimal orbits around
  the primary with a mass tending to zero.
  \citet{RobinM80} dealt with the vertical branches of planar families
  in the CRTBP and found some linearly stable \textsc{SDSP}s. 
     The accuracy of the initial values for those orbits is $10^{-6}$.
    The relation between the multiplicity and the symmetry is discussed, 
    and eight vertical branches of retrograde family $\emph{f}$ 
    around Jupiter with multicity from $5$ to $8$ are taken as examples.
    \citet{Antoniadou19} studied the stability of some
    spatial symmetric resonant periodic orbits with resonances
    $3/2, 2/1, 5/2, 3/1, 4/1, 5/1$ for both prograde and retrograde
    motions. They used the detrended Fast Lyapunov 
    Indicator (DFLI) as a tool to study the maps of dynamical stability
     (DS map) around some resonant periodic orbits.
     With similar tools, \citet{Kotoulas22} focused on the spatial symmetric 
     retrograde periodic orbits of asteroids moving in low order interior 
     mean resonances with Jupiter. 
    A review paper about the resonant periodic 
    orbits can be referred to \citet{PanSS22}. 
    
   \citet{Lara02} implemented an intrinsic predictor-corrector algorithm 
    with the help of the frenet frame. The computations of some
     periodic orbits of the CRTBP show the robust of this algorithm,
     although it is a little difficult to implement.
     \citet{Kalantonis03} transformed the solving 
    problem of the nonlinear equations into an unconstrained
    optimization problem and calculated some periodic orbits
    of Robe's CRTBP within accuracy $10^{-8}$.
    Robe's RTBP considers the motion of 
    an infinitesimal body inside one primary, which is a rigid spherical
    shell filled with a homogeneous incompressible fluid \citep{Hallan01}.
    \citet{XUAS22} independently applied Broyden's method
    with a line-search \citep{Fortran90} to the periodicity conditions 
    and implemented the scheme with the accuracy generally 
    no greater than $10^{-10}$.

    The interests on the \textsc{RTBP} come from the non-integrability
   and a lot of useful research targets, including
   three species of periodic orbits, elemental periodic orbits,
   invariant manifolds, homoclinic and heteroclinic orbits, 
   bifurcations and transit orbits. 
    For the interests on the elemental
   periodic orbits emanating from five librations and homoclinic orbits
   of the \textsc{CRTBP}, one may refer to \citet{Doedel07}. 
   The transit orbits between an interior region and an exterior region
   can be studied by the symbolic dynamics.
   Based on the values for the Sun-Jupiter-Oterma system,   
    \citet{Wilczak03} showed the existence of a homo- and heteroclinic
    cycle between two Lyapunov orbits, and also showed the existence of 
    a symbolic dynamics on four symbols. 
    \citet{Barrabes05} computed families of symmetric periodic 
    horseshoe orbits both in the planar RTBP and the spatial RTBP.
    \citet{Bengochea13} studied the numerical continuation of 
    the doubly-symmetric horseshoe orbits in the general planar
    three-body problem.
     \citet{Fitzgerald22} demonstrated the phase space geometry of
     the transit and non-transit orbits of the bicircular problem and the
     elliptic RTBP by linearing the Hamiltonian differential equations
about the collinear Lagrange points.
    For more numerical aspects on the homo- and heterclinic orbits, one may refer to \citet{Koon00,Kalantonis06,Papadakis06}
    and \citet{ZhangRY22}.

    Periodic orbits of some special versions of RTBP are studied
    recently. Sitnikov problem is a special RTBP which considers the 
    vertical motions of the infinitesimal body along a straight line
    perpendicular to the orbital plane of the primaries.    
    A brief view of the Sitnikov problem is stated in \citet{Abouelmagd20},
    and the first- and second-order approximated analytical periodic orbits 
    of the circular Sitnikov problem are constructed
    without the secular terms via the multiple scales method.
      With the background of the earth-moon-spacecraft system,
 \citet{Zaborsky20} considered the generating solutions of the problem
    of two fixed attracting masses and derived $7$ spatial families by
    continuing the angular velocity of rotation of the primaries around their 
    center of masses.
      
     Spatial Hill's problem is important for the study on the motions of
     small bodies near one primary. \citet{Zagouras85} constructed the
    fourth-order expressions for the periodic orbits originating
    at two Euler libration points near the smaller primary by
    the Lindstedt-Poincar\'{e} technique. \citet{Gomez05}
    combined the semi-analytical and numerical techniques to study the
  invariant manifold of the spatial Hill's problem associated to two Euler
  libration points.
    For a brief introduction to the recent progress on the periodic orbits
    of Hill's problem, one may refer to \citet{Kalantonis20}. 
    By computing the linear stability of the basic planar families,
    Kalantonis determined twelve vertical-self-resonant (VSR)
    periodic orbits bifurcating from planar ones,
    and found each VSR orbit generates two branches where
    the multiplicity and symmetry depend on the stability. In the
    numerical results, the periods of the spatial ones are of three or four
    times of the periods of the planar ones.   
    The quantized Hill problem can be considered as an equation system 
    derived under the effects of quantum corrections. 
     The periodic orbits emerging from the equilibrium points
     of the spatial quantized Hill problem were evaluated by the 
     averaging theory in \citet{Abouelmagd22}.  
 
   The secular behaviors of the orbits in the RTBP are important.
\citet{Prokopenya15} considered a new version of RTBP when the masses of the primaries vary isotropically with different rates, and the total mass changes according to the joint Meshcherskii law. They
studied the secular perturbations of the quasi-elliptic orbits by the 
averaging method and Hill's approximation.
   \citet{QiYi15} considered the long-term behavior of the spatial orbits
   near the Moon in the earth-moon-spacecraft CRTBP and 
   calculated some spatial periodic orbits.
   \citet{Cheng22} studied the lognormal distribution of the mass of  
   Saturn’s regular moons by the nonparameter test method in statistics, 
   and obtained the analytical expression
   of the approximate periodic orbit near the Lagrangian point $L_{3}$
   by the Lindstedt-Poincar\'{e} technique. Also, the influence some 
   parameters on the periodic orbits is also discussed.

 Although some \textsc{SDSP}s have been calculated before,
 we draw attention to the fact that our research is new, 
 including the algorithm and results.
 First, we use the numerical continuation scheme supplied in
 \citet{XUAS22}.  
 We get the desired periodicity conditions of the \textsc{SDSP}s of the autonomous
 CRTBP by the integration and Hermite interpolation, and then take use of Broyden's method with a line-search to acquire the roots as the initial values. The key program codes can be found in \citet{Fortran90}, which is
 not referenced by \citet{Kalantonis03}.
 Second, we consider the continuation of the \textsc{SDSP}s of the spatial
 Kepler problem in a rotating frame without the restriction
 on the mass ratio of the primaries, and the periods of the comet- and 
 Hill-type \textsc{SDSP}s are quite different from the known ones.
 As far as we know,
 the \textsc{SDSP}s before our work were found as vertical branches
 of some planar families or the Lyapunov orbits near the Euler equilibrium points. Third, we find $16$ initial conditions for each comet-type or Hill-type  
  resonant period ratio, and do the numerical research on the \textsc{SDSP}s systematically.

  The paper is organized as follows.
   In Sect.\ref{sec:2}, the equations of motion
   and the Hamiltonian dynamical 
   systems for the \textsc{CRTBP} and
   Hill's lunar problem are 
   introduced. In Sect.\ref{sec:3}, 
   firstly, we introduce the concept of double
   symmetry and some related lemmas.
   Secondly, the geometric way to understand 
   the comet- and Hill-type \textsc{SDSP}s of
   the approximate system is discussed.
   Thirdly, the continuation algorithm for
   the symmetric periodic orbits
   is explained. In Sect.\ref{sec:4},
   the numerical way to study the 
   linear stability of the \textsc{SDSP}s is supplied.
   In Sect.\ref{sec:5}, some numerical 
   examples are given. This article 
   ends with the discussion section.

\section{Equations and Formulation}
\label{sec:2}

   In this section, the equations of the \textsc{CRTBP} 
  for describing the comet-type orbits are introduced
  both in the inertial frame and the rotating frame, with
  the origin at the center of masses. Besides, move the origin
  to the center of one primary, the Hamiltonian for the 
  Hill-type motions is also introduced. In order to 
  get a clear relation between the rectangular coordinates
  and the canonical elliptic variables, the orbital elements,
  Delaunay elements and Poincar\'{e}-Delaunay elements 
  are introduced.
  
  \subsection{Rectangular coordinates}  
  Consider the \textsc{CRTBP} in the 
  center-of-mass inertial frame.
  Two primaries are denoted as $m_{1}$ and $m_{2}$,
  while their masses are also denoted as $m_{1}$
  and $m_{2}$, respectively. The primaries 
  move in a fixed plane which is set as the reference plane. 
  The inertial Cartesian coordinate frame
  $O-u_{1}u_{2}u_{3}$ is established by
  choosing $O$ as the center of masses and fixing one
  direction from $O$ as the $u_{1}$ axis.
  Denote the vector from  $m_{2}$ to $m_{1}$ as $\textbf{\textit{d}}$ 
  and the length $d$ is a constant. 
  The angular velocity $n^{\prime}$ of the vector 
  $\textbf{\textit{d}}$ satisfies
  $(n^{\prime})^{2}d^{3}=\mathtt{G}(m_{1}+m_{2})$,
  where $\mathtt{G}$ is the gravitational constant.
  Set the normalized units of mass, distance and time,
  $[\mathtt{M}]=m_{1}+m_{2}$, $[d]=d$ and 
  $[T]=d^{1/2}(\mathtt{G}\mathtt{M})^{-1/2}$. 
  In such units, $\mu=m_{2}/\mathtt{M}$,
  $\mathtt{G}=1$, $n^{\prime}=1$, $d=1$, and 
  $\textbf{\textit{d}}=(\cos t, \sin t, 0)^{\mathtt{T}}$,
   where the upper $\mathtt{T}$ represents transposition.  
  Let the position of the infinitesimal body be
  $\textbf{\textit{x}}=(x_{1},x_{2},x_{3}) \in\mathbb{R}^{3}
  \setminus\{\mu \textbf{\textit{d}}, (\mu-1)\textbf{\textit{d}}\}$. 
  The corresponding 
  conjugate momentum of $\textbf{\textit{x}}$ is 
  $\dot{\textbf{\textit{x}}}$.
  The differential equation system for the motion of the infinitesimal body is
  \begin{align} \label{HinerDiff}
    \ddot{\textbf{\textit{x}}}=-\frac{(1-\mu)(\textbf{\textit{x}}-\mu 
    	\textbf{\textit{d}})}{r_{1}^{3}}-\frac{\mu \left(\textbf{\textit{x}}+
    	(1-\mu)\textbf{\textit{d}}\right)}{r_{2}^{3}}.
  \end{align}  
  The Hamiltonian is composed of the kinetic energy
  and the potential energy \citep{UrsOtto}, 
  \begin{align}\label{InerHami}
   \mathcal{H}^{\texttt{iner}}=
   \frac{1}{2}\|\dot{\textbf{\textit{x}}}\|^{2}-
   \frac{1-\mu}{r_{1}}-\frac{\mu}{r_{2}} ,
  \end{align}
  where the upper ``iner" represents ``inertial", and
  $\|\cdot\|$ is the Euclidean norm,
  \begin{align*}
  & \|\dot{\textbf{\textit{x}}}\|^{2}=
  \dot{x}_{1}^{2}+\dot{x}_{2}^{2}
   +\dot{x}_{3}^{2},  \nonumber  \\
  & r_{1}^{2}=\|\textbf{\textit{x}}-\mu\textbf{\textit{d}}\|^{2}=
   \left(x_{1}-\mu \cos t \right)^{2}+
  \left(x_{2}-\mu \sin t \right)^{2}+x_{3}^{2}  ,
    \nonumber  \\
  & r_{2}^{2}=\|\textbf{\textit{x}}+(1-\mu)\textbf{\textit{d}}\|^{2}=
   \left(x_{1}+(1-\mu)\cos t \right)^{2}+
   \left(x_{2}+(1-\mu)\sin t \right)^{2}+x_{3}^{2}   .   
  \end{align*}

     Consider the \textsc{CRTBP} in 
  the center-of-mass rotating frame $O-q_{1}q_{2}q_{3}$. 
  The  vector direction of $\textbf{\textit{d}}$ is fixed 
  as the $q_{1}$ axis. 
  Denote the position vector of the infinitesimal body
  in this rotating frame as $\xi$.
  Substitute the formulas into the differential equations
  (\ref{HinerDiff}), module the rotation, 
  then the differential equations 
  in the rotating frame are acquired \citep{MeyerHO},
  \begin{align}\label{RotDiff}
   & \left(\begin{array}{c}
   \ddot{\xi}_{1}-2 \dot{\xi}_{2}   \\   
   \ddot{\xi}_{2}+2 \dot{\xi}_{1} 
   \end{array}\right)= 
   \left(\begin{array}{c}
     \xi_{1}  \\   \xi_{2}
   \end{array}\right)
   -\frac{1-\mu}{r_{1}^{3}}
   \left(\begin{array}{c}
   \xi_{1}-\mu  \\   \xi_{2}
   \end{array}\right)-\frac{\mu}{r_{2}^{3}}
   \left(\begin{array}{c}
   \xi_{1}+1-\mu  \\   \xi_{2}   
   \end{array}\right), \nonumber  \\
  & \ddot{\xi}_{3}=-\frac{\mu \xi_{3}}{r_{1}^{3}}
   -\frac{(1-\mu)\xi_{3}}{r_{2}^{3}} .
  \end{align}
    Denote $\eta_{1}=\dot{\xi}_{1}-\xi_{2}$,
  $\eta_{2}=\dot{\xi}_{2}+\xi_{1}$, and  
  $\eta_{3}=\dot{\xi}_{3}$.
  The Hamiltonian system in the rotating system can be
  checked out as follows,
  \begin{align}\label{RotHami}
  \mathcal{H}^{\texttt{rot}}=\frac{1}{2}\|\eta\|^{2}
  -n^{\prime}
  (\xi_{1}\eta_{2}-\xi_{2}\eta_{1})-\frac{1-\mu}{r_{1}}
  -\frac{\mu}{r_{2}} ,
  \end{align}  
 where $r_{1}^{2}=(\xi_{1}-\mu)^{2}+\xi_{2}^{2}+\xi_{3}^{2}$,
 and $r_{2}^{2}=(\xi_{1}+1-\mu)^{2}+\xi_{2}^{2}+\xi_{3}^{2}$.

 \subsection{Canonical elements}
  In order to better describe
  the nearly circular motions, canonical Poincar\'{e} elements are usually used. It is convenient to replace the rectangular coordinates 
  by the instantaneous orbital elements firstly, then convert the orbital elements to the Delaunay elements, and finally use 
 the Poincar\'{e}-Delaunay elements.
  
    Define two three-dimensional anticlockwise-rotating 
    matrices $\mathcal{R}_{1}$ and $\mathcal{R}_{2}$ as
   \begin{align*}
    \mathcal{R}_{1}(\theta)=\left(\begin{array}{cc}
     1  &  \text{\Large 0}      \\  
     \text{\Large 0}  &  \exp(-\textbf{\textit{J}}\theta)
    \end{array}\right), \quad
   \mathcal{R}_{2}(\theta)=\left(\begin{array}{cc}
   \exp(-\textbf{\textit{J}}\theta)  &  \text{\Large 0}   \\   %
   \text{\Large 0}  &   1
   \end{array}\right) ,
   \end{align*}
   where the skew-symmetric matrix $\textbf{\textit{J}}$ and
   the two-dimensional orthogonal matrix 
   $\exp(-\textbf{\textit{J}}\theta)$ 
   are respectively defined as 
   \begin{align*}
 \textbf{\textit{J}}=\left(\begin{array}{cc}
    0   &   1  \\
   -1   &   0   
   \end{array}\right) , \quad 
   \exp(-\textbf{\textit{J}}\theta)=\left(\begin{array}{cc}
    \cos\theta   &   -\sin\theta  \\
   \sin\theta   &   \cos\theta   
   \end{array}\right) .
   \end{align*}
   The orbital elements contains six variables, which are
   semi-major axis $a$, eccentricity $e$ ($0\leq e <1$), inclination $i$ ($0\leq i<\pi$), longitude of ascending node $\Omega$ 
   ($0\leq \Omega<2\pi$), argument of pericenter $\omega$ ($0\leq \omega<2\pi$) and mean anomaly $\ell\in \mathbb{R}$.  
   The eccentric anomaly and the true anomaly
    are $E=E(e,\ell)$ and $f=f(e,\ell)$ respectively.  
    Note that
    \begin{align*} 
     r=\|\mathbf{x}\|=\|\xi\|=a(1-e\cos E)=
     a(1-e^{2})\left(1+e\cos f\right)^{-1} .
     \end{align*} 
    Let $n$ denote the mean motion and $n$ 
    is positive for the prograde motion while negative 
    for the retrograde motion. For the unperturbed system,
    if the mass of the central primary is $\mu$, then 
    we have $n^{2}a^{3}=\mu$.
    The angular momentum is denoted 
    as $\Theta=\pm \sqrt{ma(1-e^{2})}$, 
    where the ``$+$" represents the 
    prograde motion while ``$-$" represents 
    the retrograde motion.     
    The rectangular coordinates of the position
    can be expressed by
    \begin{align*}
     \mathbf{x}=\mathcal{R}_{2}(\Omega)
     \mathcal{R}_{1}(i) 
     \mathcal{R}_{2}(f+\omega)
     \left(r,0,0\right)^{\mathtt{T}} .
    \end{align*}
    In the sense of instantaneous, the velocity is about
    the derivative of the anomalies,
    \begin{align*}
     \dot{\mathbf{x}}=\mathcal{R}_{2}(\Omega)
     \mathcal{R}_{1}(i) \mathcal{R}_{2}(f+\omega)
     \left(\dot{r},r\dot{f},0\right)^{\mathtt{T}},
    \end{align*}
   where $\dot{r}= e\sin(f)/\Theta$, 
    $r\dot{f}= (1+e\cos f)/\Theta $. 
   According to the formulas below, the anomalies can
   be converted to each other,
   \begin{align*}
    E-e\sin E=\ell, \quad r\cos f=a(\cos E -e), \quad 
    r\sin f=a\sqrt{1-e^{2}}\sin E .
   \end{align*}    
   The Delaunay elements are
   \begin{eqnarray}
   \begin{array}{ccccccccc}
   L    &=& \sqrt{ \mu a}, &  G &=&  L\sqrt{1-e^{2}},  &  H &=&  G\cos i,  \nonumber  \\
   \ell &=& n t+\ell_{0},   &  g &=&  \omega,    &  h &=&  \Omega,
   \end{array}
   \end{eqnarray}
   and the Poincar\'{e}-Delaunay elements are
   \begin{align}
   \begin{array}{ccccccccc}
    P_{1} &=&  L-G+H,    &  P_{2}  &=&  \sqrt{2(L-G)}\cos(g+h),  &  P_{3} &=&   G-H ,
    \nonumber  \\
   Q_{1} &=&  \ell+g+h, &  Q_{2}  &=&  -\sqrt{2(L-G)}\sin(g+h), &  Q_{3} &=&  \ell+g .  
   \end{array}
   \end{align}
   
   \subsection{Perturbed system}
   Consider the comet-type motion of the infinitesimal body,
   so $1/r$ is small.
   Let $\{\mathcal{P}_{k}\}_{k=0}^{\infty}$ be the 
   Legendre polynomials, and $\mathcal{P}_{0}=1$,
   $\mathcal{P}_{1}(\xi_{1})=\xi_{1}$,
   $\mathcal{P}_{2}(\xi_{1})=\frac{3}{2}\xi_{1}^{2}
   -\frac{1}{2}$, $\mathcal{P}_{3}(\xi_{1})=
   \frac{5}{2}\xi_{1}^{3}-\frac{3}{2}\xi_{1}$.
   According to the Legendre expansion technique,   
   \begin{align*}
    \frac{1-\mu}{r_{1}}+\frac{\mu}{r_{2}} &=
    \frac{1-\mu}{\sqrt{r^{2}+\mu^{2}-2\mu \xi_{1}}}
    +\frac{\mu}{\sqrt{r^{2}+(1-\mu)^{2}+2(1-\mu) \xi_{1}}}
    \nonumber  \\
    &=\frac{1}{r}+\frac{1-\mu}{r}\sum_{j=2}^{\infty}
    \mathcal{P}_{j}\left(\frac{\xi_{1}}{r}\right)
    \left(\frac{\mu}{r}\right)^{j}+
    \frac{\mu}{r}\sum_{k=2}^{\infty}
    \mathcal{P}_{k}\left(-\frac{\xi_{1}}{r}\right)
    \left(\frac{1-\mu}{r}\right)^{k} .
   \end{align*}
   If $1/r$ is sufficiently small, the Hamiltonian 
   $\mathcal{H}^{\texttt{rot}}$ can be
   written in the perturbed form as
   $\mathcal{H}^{\texttt{rot}}=\mathcal{H}_{0}^{\texttt{rot}}
   +\mathcal{O}(1/r^{3})$ by the canonical elements, where
   \begin{align}\label{H0rot}
    \mathcal{H}_{0}^{\texttt{rot}}=
   -\frac{1}{2L^{2}}-H=-\frac{1}{2(P_{1}+P_{3})^{2}}
   -\left(P_{1}-\frac{P_{2}^{2}+Q_{2}^{2}}{2}\right) .   
   \end{align}

    If the infinitesimal body is very close to primary $m_{2}$,
   it is convenient to move the origin to the center of $m_{2}$. 
   Let $m_{2}-q_{1}q_{2}q_{3}$ denote
   the $m_{2}$ centered rotating frame.
   Let the position $\textbf{\textit{q}}=\xi+(1-\mu,0,0)^{\mathtt{T}}$, 
   the conjugate momentum $\textbf{\textit{p}}=\eta+
   (0,1-\mu,0)^{\mathtt{T}}$. Note that $\|\textbf{\textit{q}}\| \ll 1$. 
   The Hamiltonian for the Hill-type 
   motion around $m_{2}$ can be written as
   \begin{align*}
   \mathcal{H}^{\texttt{mtwo}}&=\mathcal{H}_{0}^{\mu}
   +(1-\mu) q_{1}-\frac{1-\mu}{\sqrt{1-2q_{1}
   +\|\textbf{\textit{q}}\|^{2}}}
   -\frac{(1-\mu)^{2}}{2} ,  \nonumber  \\
    \mathcal{H}_{0}^{\mu}&=\frac{1}{2}\|\textbf{\textit{p}}\|^{2}-
    (q_{1}p_{2}- q_{2}p_{1})-\frac{\mu}{\|\textbf{\textit{q}}\|} .
   \end{align*}
   Expand the potential function by Legendre expansion
   technique, neglect the constant $-(1-\mu)-(1-\mu)^{2}/2$,
   and write the part of the integrable Hamiltonian in the canonical elements,
    we get
   \begin{align}\label{Hmtwo}
   \mathcal{H}^{\texttt{mtwo}}=\mathcal{H}_{0}^{\mu}
   -(1-\mu) \left(\frac{3}{2}q_{1}^{2}-\frac{1}{2}
   \|\textbf{\textit{q}}\|^{2}\right)
   -(1-\mu) \cdot \mathcal{O}(\|\textbf{\textit{q}}\|^{3}) ,
   \end{align}   
   where $\mathcal{H}_{0}^{\mu}=-\frac{\mu}{2L^{2}}-H$ 
   and the canonical elements change with $ L= \sqrt{ \mu a}$.
    
    Hill's lunar problem can be derived by the following way.  
   Make symplectic scaling 
   $\textbf{\textit{q}}\rightarrow \mu^{1/3}\textbf{\textit{q}}$
   and $\textbf{\textit{p}}\rightarrow \mu^{1/3}\textbf{\textit{p}}$, 
   then scale the Hamiltonian by $\mu^{-2/3}$, afterwards make 
   $\mu\rightarrow 0$, the Hamiltonian becomes
   \begin{align}\label{HaHill}
   \mathcal{H}^{\texttt{Hill}}=\mathcal{H}_{0}^{\texttt{rot}}
     -\left(\frac{3}{2}q_{1}^{2}-\frac{1}{2}
     \|\textbf{\textit{q}}\|^{2}\right) .
   \end{align}
  In the following, the integrable system like
  $\mathcal{H}_{0}^{\texttt{rot}}$ and $\mathcal{H}_{0}^{\mu}$
  is usually referred as the \emph{approximate system}.
 
\section{Double symmetries and \textsc{SDSP}s}
\label{sec:3}
  
    The double symmetries exist in the \textsc{RTBP} and 
   Hill's lunar problem, because these systems
   satisfy two time-reversing symmetries. One is 
   with respect to the syzygy (the line containing 
   both primaries), the other
   is about the vertical plane containing the primaries.     
  The differential equation system (\ref{RotDiff}) 
  is invariant under the two anti-symplectic reflections
  \citep{Kazantis79,RobinM80,Howison00I,Howison00II,
  	LlibreR09,LlibreS11,Bengochea13,XbXu19,XuHill20,
 Antoniadou19,Kotoulas22}.
  \begin{align*}
 \mathscr{R}_{1}: & (\xi_{1},\xi_{2},\xi_{3},
 \dot{\xi}_{1}- \xi_{2},
 \dot{\xi}_{2}+ \xi_{1},
 \dot{\xi}_{3})
 \rightarrow  (\xi_{1},-\xi_{2},-\xi_{3},
 -\dot{\xi}_{1}+ \xi_{2},
 \dot{\xi}_{2}+ \xi_{1},
 \dot{\xi}_{3}), \nonumber  \\
 \mathscr{R}_{2}:  &
 (\xi_{1},\xi_{2},\xi_{3},
 \dot{\xi}_{1}- \xi_{2},
 \dot{\xi}_{2}+ \xi_{1},
 \dot{\xi}_{3})
 \rightarrow         (\xi_{1},-\xi_{2},\xi_{3},
 -\dot{\xi}_{1}+ \xi_{2},
 \dot{\xi}_{2}+ \xi_{1},
 -\dot{\xi}_{3}).
 \end{align*}
 The two time reversing symmetries can be explained 
 geometrically. $\mathscr{R}_{1}$ denotes the 
 symmetry about the $\xi_{1}$ axis, while $\mathscr{R}_{2}$
 denotes the symmetry about the $\xi_{1}\xi_{3}$ plane. 
  It is difficult to know directly the rectangular initial values for the \textsc{SDSP}s, but there is one way to
 get the precise initial values by numerically continuing
 the approximate initial values, which can be achieved by determining the \textsc{SDSP}s of the approximate system.
 
  Though the Hamiltonians of the comet and Hill-type
 motions are different, the corresponding approximate
 systems are similar to $\mathcal{H}_{0}^{\texttt{rot}}$. 
 The \textsc{SDSP}s of the approximate system can be
 determined according to the following three lemmas
 \citep{Howison00I,XbXu19}.
 \begin{lemma}
 If one orbit of the Hamiltonian (\ref{RotHami}) starts from a Lagrangian
 subplane $\mathscr{L}_{1}^{(0)}$ (or $\mathscr{L}_{2}^{(0)}$),
 and intersects $\mathscr{L}_{2}^{(0)}$ (or $\mathscr{L}_{1}^{(0)}$) after time $T/4>0$, then the orbit is periodic with
 period $T$ and doubly-symmetric, where
 \begin{align}\label{boundaryVal}
 \mathscr{L}_{1}^{(0)} &=\{(\xi,\eta)|\xi=(\xi_{1},0,0),
 \eta=(0,\dot{\xi}_{2}+\xi_{1},\dot{\xi}_{3})\},  
 \nonumber \\
 \mathscr{L}_{2}^{(0)} &=\{(\xi,\eta)|\xi=(\xi_{1},0,\xi_{3}),
 \eta=(0,\dot{\xi}_{2}+\xi_{1},0) \} .
 \end{align}
 The Hamiltonian systems (\ref{H0rot}), (\ref{Hmtwo}), 
 (\ref{HaHill}) also satisfy this proposition.
 \end{lemma}  
\begin{proof}
  Consider one solution starts from $\mathscr{L}_{1}^{(0)}$.
  Let $\phi(t,X_{0},Y_{0})$ be a solution of 
  (\ref{RotDiff}) with the initial values 
  $X_{0}=(\xi_{1},0,0)$ and	$Y_{0}=(0,\eta_{2},\eta_{3})$ 
  at the time $t=0$. 
  The differential equation system satisfies the 
  double symmetries $\mathscr{R}_{1}$ and $\mathscr{R}_{2}$.
  After time $T/4$, the solution 
  intersects with the set $\mathscr{L}_{2}^{(0)}$
  at $X_{1}=(\tilde{\xi}_{1},0,
  \tilde{\xi}_{3})$,
  $Y_{1}=(0,\tilde{\eta}_{2},0)$.
  One has 
  \begin{align*}
   \phi_{T/4}=\phi(T/4,X_{0},Y_{0})=(\tilde{\xi}_{1},0,
  \tilde{\xi}_{3},
  0,\tilde{\eta}_{2},0) .
  \end{align*}
  As the solution satisfies the symmetry $\mathscr{R}_{1}$,
  one has
  \begin{align*}
    \mathscr{R}_{1} \phi_{T/4}
   =\phi(-T/4,X_{0},Y_{0})
   =(\tilde{\xi}_{1},0,-\tilde{\xi}_{3},
  0,\tilde{\eta}_{2},0) .	
  \end{align*} 
  It takes $T/2$ for this orbit from $\mathscr{L}_{2}$ to 
  $\mathscr{L}_{2}$ again.
  The solution also satisfies the symmetry $\mathscr{R}_{2}$, so
  one has
  \begin{align*}
  \phi(T/4+T/2,X_{0},Y_{0})=\mathscr{R}_{2}\phi_{-T/4}
  =(\tilde{\xi}_{1},0,-\tilde{\xi}_{3},
  0,\tilde{\eta}_{2},0)
  =\phi(-T/4,X_{0},Y_{0}) .
  \end{align*}
  So the solution $\phi(t,X_{0},Y_{0})$ is periodic with
  period $T$ and doubly symmetric.  
  For the case that one solution starts
  from $\mathscr{L}_{2}^{(0)}$, the proof is the same as above.
\end{proof}
\begin{lemma}\label{Lem2}
  The Lagrangian subplane $\mathscr{L}_{k}^{(0)}$ is equivalent
  to $\mathscr{L}_{k}^{(1)}$, to $\mathscr{L}_{k}^{(2)}$
  and to $\mathscr{L}_{k}^{(3)}$ for $k=1, 2$, where
   \begin{align}
   \mathscr{L}_{1}^{(1)} &=\{(f=0 \mod \pi, \ \ \omega=0 \mod \pi, \ \ \Omega=0 \mod \pi )\},  \nonumber  \\
   \mathscr{L}_{2}^{(1)} &=\{(f=0 \mod \pi, \ \ \omega=\frac{\pi}{2} \mod \pi, \ \ \Omega=\frac{\pi}{2} \mod \pi ) \}.
   \end{align}
   \begin{align}
   \mathscr{L}_{1}^{(2)} &=\{(\ell=0 \mod \pi, \ \  g=0 \mod \pi, \ \                h=0 \mod \pi )\},  \nonumber  \\
   \mathscr{L}_{2}^{(2)} &=\{(\ell=0 \mod \pi, \ \  g=\frac{\pi}{2} \mod \pi, \ \    h=\frac{\pi}{2} \mod \pi ) \}.
   \end{align}
   \begin{align}\label{PoinSymmetry}
    \mathscr{L}_{1}^{(3)} &=\{ (Q_{1}=0 \mod \pi,  \  \  Q_{2}\equiv 0,  \ \ Q_{3}=0 \mod \pi  )\},  \nonumber \\
    \mathscr{L}_{2}^{(3)} &=\{ (Q_{1}=0 \mod \pi,  \  \  Q_{2}\equiv 0,  \ \ Q_{3}=\frac{\pi}{2} \mod \pi )\}.
   \end{align}
\end{lemma}
 \begin{lemma}\label{Lem3}
 	The \textsc{SDSP}s of the approximate system 
 	$\mathcal{H}_{0}^{\texttt{rot}}$ can only be the circular orbits.
\end{lemma}
\begin{proof}
The rectangular coordinates $\xi$ and $\eta$ can be expressed
by the orbital elements,
\begin{align*}
 \xi =\mathcal{R}_{2}(-n^{\prime}t)\textbf{\textit{x}}
 =r\alpha_{1}, \quad  
 \eta =\mathcal{R}_{2}(-n^{\prime}t)\dot{\textbf{\textit{x}}}
 =\dot{r}\alpha_{1}+r\dot{f}\alpha_{2},  
\end{align*}
 where $\mathcal{A}=\mathcal{R}_{2}
 (\Omega-n^{\prime}t)
 \mathcal{R}_{1}(i) \mathcal{R}_{2}(f+\omega)
 =(\alpha_{1},\alpha_{2},\alpha_{3})$ 
 is a $3\times 3$ matrix,
 \begin{align*}
\alpha_{1}&=\left(\begin{array}{c}
	\cos(f+\omega)\cos(\Omega-t)-\cos i \sin(f+\omega)
	\sin(\Omega-t)     \\
	\cos(f+\omega)\sin(\Omega-t)+\cos i \sin(f+\omega)
	\cos(\Omega-t)     \\
	\sin(f+\omega)\sin i
\end{array}\right),  \\
\alpha_{2}&=\left(\begin{array}{c}
-\sin(f+\omega)\cos(\Omega-t)-\cos i \cos(f+\omega)
\sin(\Omega-t)     \\
-\sin(f+\omega)\sin(\Omega-t)+\cos i \cos(f+\omega)
\cos(\Omega-t)     \\
\cos(f+\omega)\sin i
\end{array}\right),  
\end{align*}
and $\alpha_{3}=\left(\sin(\Omega-t)\sin i, 
-\cos(\Omega-t)\sin i, \cos i \right)^{\mathtt{T}}$.
 Let the initial values of one orbit be
 $X_{0},Y_{0}$ as above. Suppose that after time $T/4$,
 the solution intersect $\mathscr{L}_{2}^{(0)}$ at
 $X_{1},Y_{1}$. 
 According to the two boundary conditions, 
 the solution satisfies 
 $f=\omega=\Omega=0 \mod \pi$ at $t=0$, and
 $f=0 \mod \pi$, $g=\pi/2 \mod \pi$, $h=\pi/2 \mod \pi$
 at $t=T/4$. However, $\dot{g}=0$ in the approximate system.
 There exists a conflict if the orbit is elliptic.
 So the doubly-symmetric periodic orbit can only be circular.
\end{proof}

  \section{Approximate system and \textsc{SDSP}s}
  \label{sec:4}
    Let us recall three theorems proved by \citet{Howison00I,Howison00II}
 in order to make the readers to understand what kinds of 
 \textsc{SDSP}s to be calculated in this paper. 
 \begin{theorem}{\citep{Howison00I,Howison00II}}
  (1) There exist doubly-symmetric periodic solutions of the spatial
restricted three-body problems for all values of the mass ratio 
parameter $\mu$ with large inclination which are arbitrarily far away from the primaries.
 (2) There exist doubly-symmetric periodic solutions of the spatial
restricted three-body problems for all values of the mass ratio parameter $\mu$ with large inclination which are arbitrarily close to one of the primaries.
 (3) There exist doubly-symmetric periodic solutions of the spatial Hill’s lunar problems with large inclination which are arbitrarily close to the primary.
   \end{theorem}
   The existence of the \textsc{SDSP}s is shown by
   the continuation method, which needs a small parameter.
   The small parameter is introduced by the symplectic
   scaling method\citep{Howison00I}. 
   Theoretically, the small parameter is small enough, however,
   it is not necessarily to be too small in the numerical research.
   
   \subsection{Symplectic scaling}
   In order to understand the small parameter and estimate
   the order of magnitude of the perturbation to the approximate
   system, we recall the symplectic scaling procedure\citep{MeyerHO,
   Howison00I,XbXu19,XuHill20}.   
   For the comet-type orbits, the 
   coordinate frame is the rotating frame with the 
   origin at the center of masses. The symplectic transformation is
   $ \xi \rightarrow \varepsilon^{-2}\xi, \quad 
  \eta \rightarrow \varepsilon \eta $,
   and the new Hamiltonian is equal to the old Hamiltonian
  multiplied with $\varepsilon$. The small parameter $\varepsilon^{2}$ represents the inverse of the great distance of the infinitesimal body 
   from the origin $1/r$. 
   The symplectic scaled new Hamiltonian 
   is $\tilde{\mathcal{H}}^{\texttt{rot}}
  =\tilde{\mathcal{H}}_{0}^{\texttt{rot}}
  +\mathcal{O}(\varepsilon^{7})$,
  where $\tilde{\mathcal{H}}_{0}^{\texttt{rot}}=
  -\varepsilon^{3}/(2\tilde{L}^{2})-\tilde{H}$, $L=\sqrt{a}$, 
  $\tilde{L}=\varepsilon L$ and $\tilde{H}=\varepsilon H$.
  
    While for the Hill-type orbits, the origin is set to locate at 
   the center of one primary $m_{2}$. This symplectic transfromation is
    $ \xi \rightarrow \varepsilon^{2}\xi, \quad 
    \eta \rightarrow \varepsilon^{-1}\eta $,
    and the new Hamiltonian is equal to 
    the old Hamiltonian multiplied with $\varepsilon^{-1}$.
     Here the small parameter $\varepsilon^{2}$ represents the degree of 
   closeness to the primary $m_{2}$. 
   $\mathcal{H}^{\texttt{mtwo}}$ is transformed to
   $\hat{\mathcal{H}}^{\texttt{mtwo}}=
   \hat{\mathcal{H}}_{0}^{\mu}+\mathcal{O}(\varepsilon^{3})$, 
   where $ \hat{\mathcal{H}}_{0}^{\mu}
   =-\frac{\varepsilon^{-3}\mu}{2\hat{L}^{2}}-\hat{H}$, 
   $L=\sqrt{\mu a}$, $\hat{L}=\varepsilon^{-1} L$ and $\hat{H}
   =\hat{L}\sqrt{1-e^{2}}\cos i$. If $\mu\rightarrow 0$,
   we use $\hat{L}= \mu^{1/3}\check{L}$, and have
   $\check{\mathcal{H}}^{\texttt{Hill}}=
   \check{\mathcal{H}}_{0}^{\texttt{rot}}+\mathcal{O}
   (\varepsilon^{3})$, where $\check{\mathcal{H}}_{0}^{\texttt{rot}}
   = -\varepsilon^{-3}/(2\check{L}^{2})-\check{H}$.
     
   \subsection{Sixteen cases of the initial values}
    In order to do the numerical study, it is convenient to use the usual rectangular coordinates to represent the initial values. 
   Consider the initial values of the comet-type \textsc{SDSP}s of the approximate system. Let $\phi(t,X_{0},Y_{0})$ be one solution and satisfy
   $(X_{0},Y_{0}) \in \mathscr{L}_{1}^{(0)}$.
   According to Lemma \ref{Lem2} and Lemma \ref{Lem3},
   the initial values at $t=0$ can be written as
  \begin{align*}
   	X_{0}=(\pm a\cos\Omega_{0},0,0), 
   	\quad
   Y_{0}
   =\pm (0, n a\cos i\cos\Omega_{0}, 
   n a\sin i) . 
  \end{align*}
  Suppose $\phi(\frac{T}{4},X_{0},Y_{0})
  =(X_{1},Y_{1})$$\in\mathscr{L}_{2}^{(0)}$. Let 
  $\Omega_{1}=\Omega_{0}-T/4$, and 
  we have again
   \begin{align*}
 X_{1}=a\left(\mp\cos i\sin(\Omega_{1}), 0, \pm\sin i \right), \quad 
 Y_{1}=\left(0, \pm n a \sin(\Omega_{1}), 0\right). 
  \end{align*}
  Usually, $\dot{X}_{k}=Y_{k}-(0,a,0)$  for $k=1, 2$ 
  are used instead of $Y_{0}, Y_{1}$, respectively.
  
\begin{figure}[h]
   	\center{\includegraphics [scale=0.47]{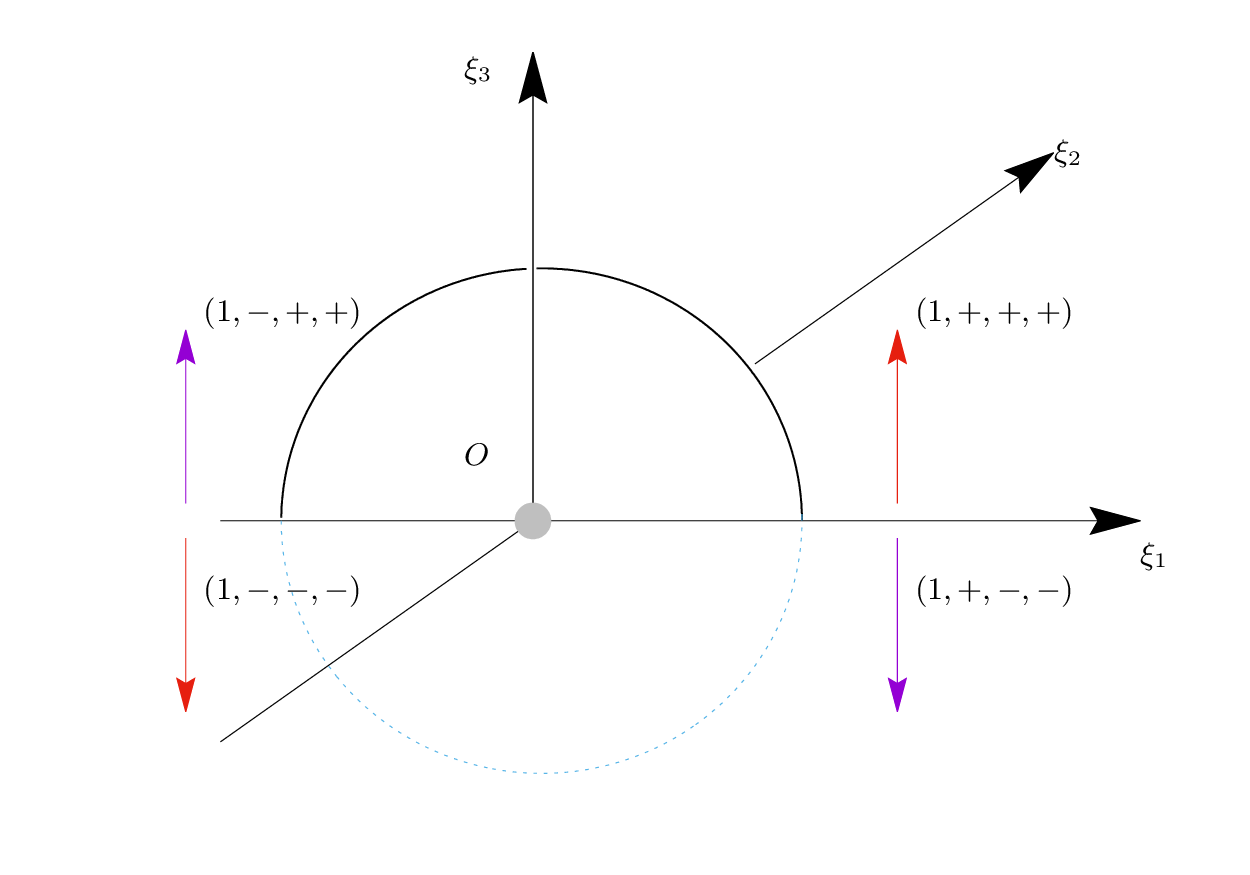}
   	\includegraphics [scale=0.47]{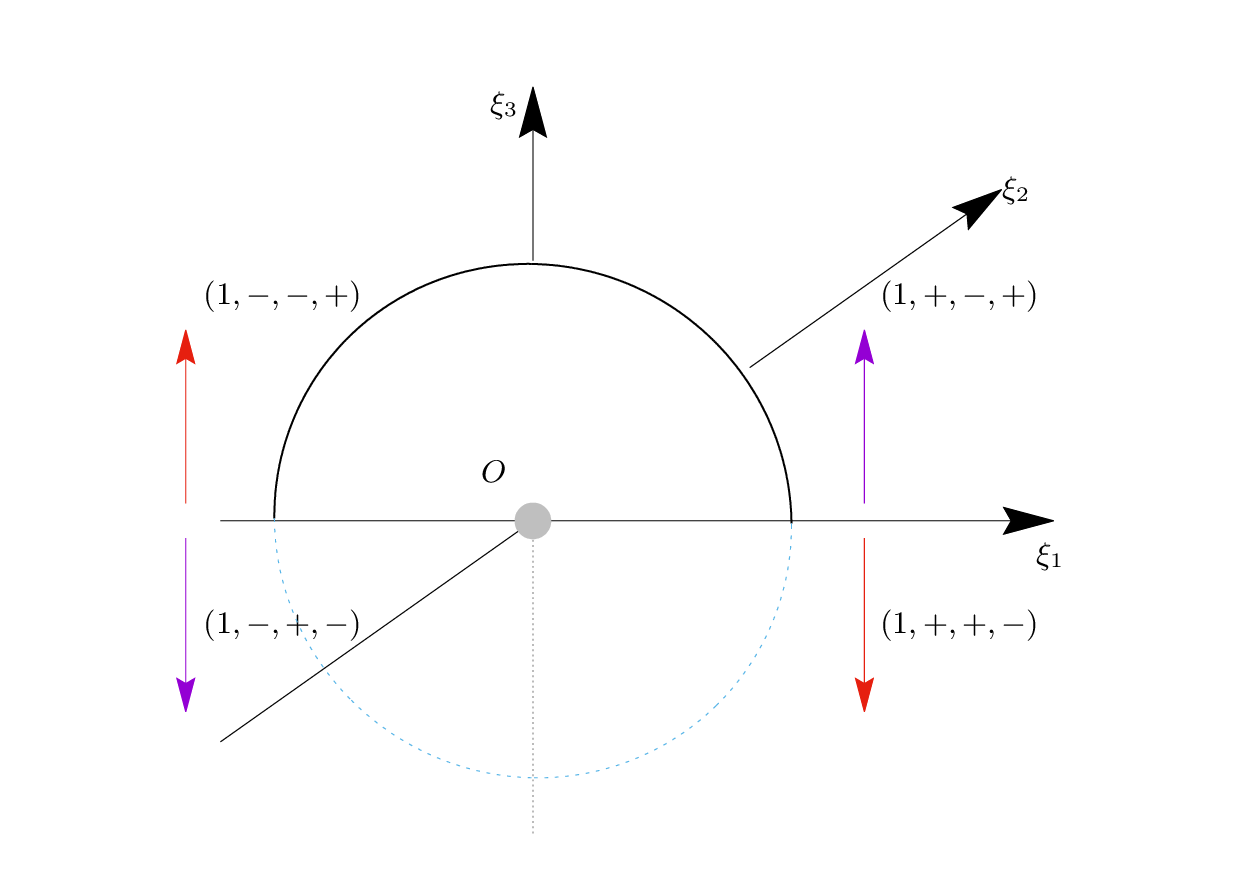}
   	 	\includegraphics [scale=0.47]{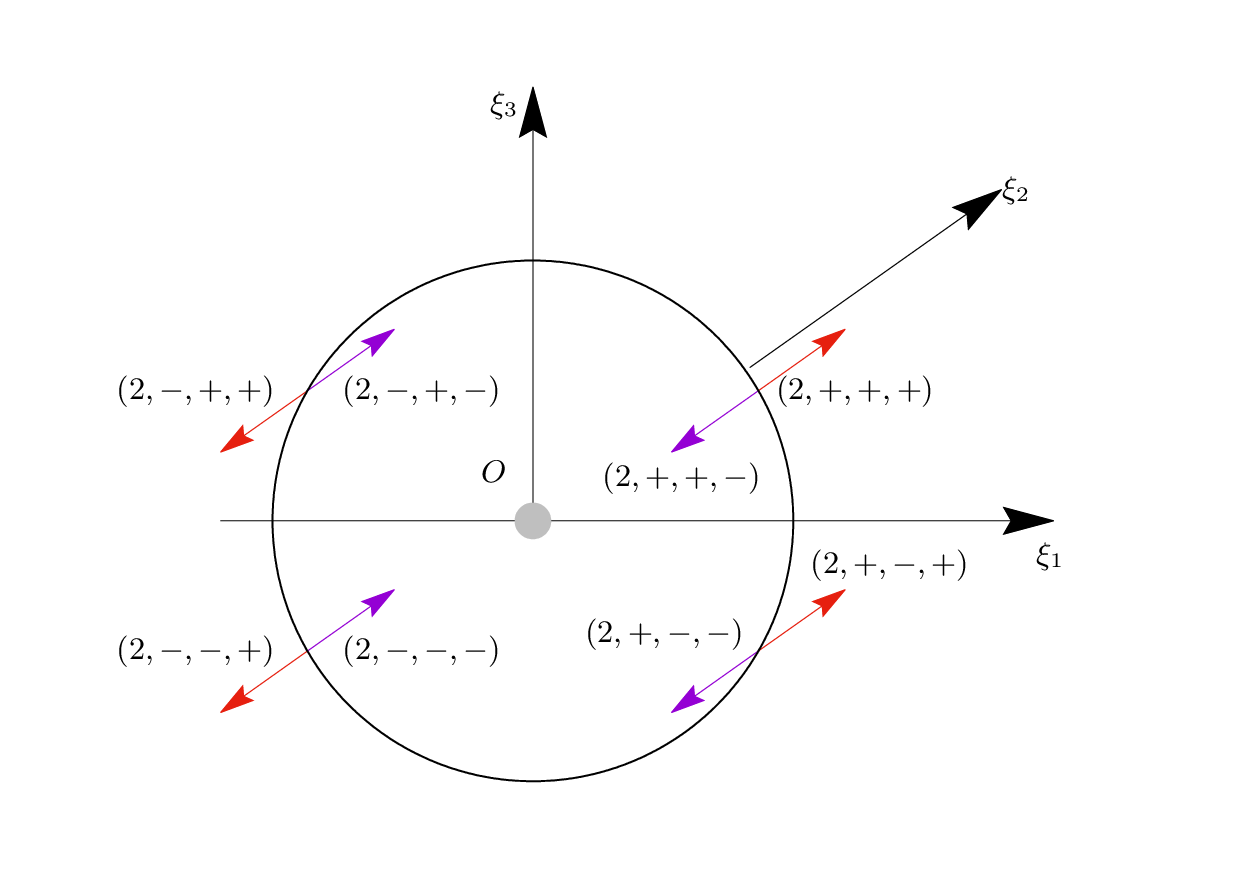}}
   	\caption{ Three sketch maps are supplied for describing the initial values for the \textsc{SDSP}s of the approximate system. 
   	The included angle between the initial orbital plane 
   	$\mathscr{L}_{1}^{(0)}$ and the $\xi_{2}$ axis is no greater
   	than $90^{\circ}$ for the first sketch (upper left) and 
   	greater than  $90^{\circ}$ for the second sketch (upper right). 
   	There are four cases of initial values in each of the 
   	first two sketch maps. In the third sketch map (lower middle), 
   	the included angle between
   	the initial orbital plane $\mathscr{L}_{2}^{(0)}$ and the 
   	$\xi_{2}$ axis is $90^{\circ}$, and there are also eight cases of 
   	initial values. The signs of $(\xi_{1},\dot{\xi}_{2}, \dot{\xi}_{3})$ 
   	are marked in the upper two sketch maps, while the signs of $(\tilde{\xi}_{1},\tilde{\xi}_{3},
   	\dot{\tilde{\xi}}_{2})$ are marked in the
   	third sketch map. There are totally
  	sixteen cases of initial values for the \textsc{SDSP}s 
  	of the approximate system. }\label{fig1IV2}
   \end{figure}

  Let $\mathrm{sgn}(n)=1$ if the infinitesimal body 
  moves prograde (anticlockwise) and $\mathrm{sgn}(n)=-1$ 
  if the infinitesimal body moves retrograde (clockwise).
  Given an acute angle $i_{0}<\pi/2$,
  the inclination is supposed to be $i_{0}$.
  If the orbit starts from $\mathscr{L}_{1}^{(0)}$,
  each parameter in $(\xi_{1},\dot{\xi}_{2},\dot{\xi}_{3})$ 
  has a possibility to be positive or negative, 
  so there are eight cases of the initial values.
  There are also eight cases for the signs of $(\xi_{1},
  \xi_{3},\dot{\xi}_{2})$ if the orbit starts from 
  $\mathscr{L}_{2}^{(0)}$.
  As we can see in Fig. \ref{fig1IV2},
  there are three orbital planes with different dihedral angles.
   The first two sketch maps represent that the orbits start 
   from $\mathscr{L}_{1}^{(0)}$, and the third sketch map represents
   that the orbits start from $\mathscr{L}_{2}^{(0)}$. There are eight cases of initial values in each Lagrangian subplane. 
   With the orbital elements, we use 
   The signs of $(\xi_{1},\dot{\xi}_{2},\dot{\xi}_{3})$ or $(\xi_{1},\xi_{3},\dot{xi}_{2})$ are applied to represent a set of initial values. 
   For example, $(1,+,+,+)$
   means that the orbit starts from 
   $\mathscr{L}_{1}^{(0)}$ and
   $\xi_{1}=a>0$, $\dot{\xi}_{2}>0$, 
   $\dot{\xi}_{3}>0$, while $(2,+,+,+)$ represents
   that the orbit starts from $\mathscr{L}_{2}^{(0)}$
   and $\tilde{\xi}_{1}=a\cos i_{0}>0$, 
   $\tilde{\xi}_{3}>0$,$\dot{\tilde{\xi}}_{2}>0$. 
   In all, there can be $16$
   possibilities of initial values of \textsc{SDSP}s 
   for a fixed period ratio in the approximate system. 
   For more details, the readers can refer 
   to the caption of the Fig. \ref{fig1IV2}.

  \subsection{Initial values}
   According to the Hamiltonian $\mathcal{H}_{0}^{\texttt{rot}}$
  of the approximate system in Eq.(\ref{H0rot}), the 
  longitude of the ascending node satisfies the
  differential equation $
    \dot{h}=\frac{\partial \mathcal{H}_{0}^{\texttt{rot}}}
    {\partial H}=-1 $, 
  so $h(t)=h(0)-t$ and the line of apsides moves retrograde.  
  Let $\mathbb{Z}^{+}$ denote all the non-negative integers.
  One fourth period is set to be $\frac{T_{0}}{4}=\frac{\pi}{2}+k\pi$ 
  ($k\in \mathbb{Z}^{+}$), such that 
   $h(\frac{T_{0}}{4})=h(0)-\frac{T_{0}}{4}$ satisfies 
  $\mathscr{L}_{2}^{(2)}$ if $h(0)=0 \mod \pi$ and
  satisfies $\mathscr{L}_{1}^{(2)}$ if $h(0)=\frac{\pi}{2} \mod \pi$.
  In addition,  we have the differential equation
  $\dot{Q}_{3}(t)=\dot{\ell}(t)+\dot{g}(t)=
  \mathrm{sgn}(n) L^{-3}$, so 
  $Q_{3}(t)=Q_{3}(0)\pm L^{-3}t$. We must have
  $L^{-3}\frac{T_{0}}{4}=\frac{\pi}{2}+j\pi $ with 
  $j\in \mathbb{Z}^{+}$, such that 
  $Q_{3}(\frac{T_{0}}{4})$ satisfies $\mathscr{L}_{2}^{(3)}$
  if $Q_{3}(0)=0 \mod \pi$ and satisfies $\mathscr{L}_{1}^{(3)}$
  if $Q_{3}(0)=\frac{\pi}{2} \mod \pi$. Because  
   the argument of perigee $g$ is of no definition in the circular orbit,
  and we redefine $\ell(\frac{T_{0}}{4})
  =\ell(0)+\mathrm{sgn}(n)j\pi$ and 
  $g(\frac{T_{0}}{4})=\mathrm{sgn}(n)\frac{\pi}{2}$
  such that this definition keeps in coincident with
  the Poincar\'{e}-Delaunay elements for the symmetries.   
  Then $T_{0}=2\pi+4k\pi$ is set as the period of a \textsc{SDSP} of the 
  approximate system. In the inertial frame,  
  during the period that the primaries revolve $2k+1$ turns around each 
  other, the projection of infinitesimal body in the $\xi_{1}\xi_{2}$ plane  
  revolve $2j+1=\frac{T_{0}}{2\pi}|n|$ turns.

   The mean anomaly can be determined by the period ratio, that is
   $|n|=\frac{2j+1}{2k+1}$.
  For the comet-type orbits, the semi-major axis $a$ is 
  much longer than $d=1$ and is determined by
  $a=L^{2}=\left[(2k+1)/(2j+1)\right]^{2/3}$ with $k\gg j\geq 0$. 
  We can check
  that the Kepler's third law $n^{2}a^{3}=1$ is satisfied for the 
  comet-type \textsc{SDSP}s.  
   The small parameter $\varepsilon$ is also determined by
   the period ratio as
    $\varepsilon^{3}=\left[(2j+1)/(2k+1)\right]$ with $\tilde{L}=1$.    
    We conclude that $n a=\mathrm{sgn}(n) \varepsilon
    =a^{-1/2}$. Suppose $a=a_{0}$ when $j$ and $k$ are fixed.     
    The initial orbital inclination is $i_{0}<\frac{\pi}{2}$ or $\pi-i_{0}$. 
    Then the initial values of a \textsc{SDSP} starting from
    $\mathscr{L}_{1}^{(0)}$ can be written as
    \begin{align}\label{InitV}
    (\xi_{1},0,0,0,\dot{\xi}_{2},\dot{\xi}_{3})=
    (\pm a_{0},0,0,0,\pm a_{0}^{-1/2}\cos i_{0} -a_{0}, 
    \pm a_{0}^{-1/2}\sin i_{0}). 
  \end{align} 
  While the initial values of a \textsc{SDSP} starting from 
  $\mathscr{L}_{2}^{(0)}$ can be written as
  \begin{align}\label{InitV2}
    (\tilde{\xi}_{1},0,\tilde{\xi}_{3},0,\dot{\tilde{\xi}}_{2},0)=
     (\pm a_{0}\cos i_{0}, 0, \pm a_{0}\sin i_{0}, 0,
    \pm a_{0}^{-1/2} -a_{0}, 0).  
  \end{align}

  The initial values $(\mathbf{q},\dot{\mathbf{q}})$ of a Hill-type 
  \textsc{SDSP} of $\mathcal{H}_{0}^{\mu}$ in (\ref{Hmtwo}) are 
  formally the same as (\ref{InitV}) or (\ref{InitV2}). 
  For the Hill-type orbits, 
  the semi-major axis $a_{1}$ is much shorter than $d=1$,
  and $a_{1}=\left(\frac{\mu}{n^{2}}\right)^{\frac{1}{3}}
  =\mu^{1/3}\left[(2k+1)/(2j+1)\right]^{2/3}$ with $j \gg k \geq 0$. 
  Denote the initial orbital inclination as $i_{1}<\frac{\pi}{2}$ 
  or $\pi-i_{1}$.  
  The initial parameters of a Hill-type \textsc{SDSP} starting from
  $\mathscr{L}_{1}^{(0)}$ can be expressed as
  \begin{align} \label{Hill14}
  (q_{1},\dot{q}_{2},\dot{q}_{3})=
  (\pm a_{1}, \pm\sqrt{\mu/a_{1}}\cos i_{1} -a_{1},
  \pm \sqrt{\mu/a_{1}}\sin i_{1}) .
  \end{align}
  If a Hill-type \textsc{SDSP} starts from $\mathscr{L}_{2}^{(0)}$,
  the set of initial parameters is given by
   \begin{align} \label{Hill15}
  (\tilde{q}_{1},\tilde{q}_{3},\dot{\tilde{q}}_{2})=
  (\pm a_{1}\cos i_{1}, \pm a_{1}\sin i_{1},
  \pm \sqrt{\mu/a_{1}}-a_{1}) .
  \end{align}  
  The initial values for the Hill-type \textsc{SDSP}s of the 
  approximate system of Hill's lunar problem 
  can be calculated by Eqs.(\ref{InitV}) and (\ref{InitV2}) 
  for $j \gg k \geq 0$.   
  
  \subsection{Periodicity conditions}
    The integration procedure and the periodicity conditions 
    are based on the rectangular coordinates. 
    Consider the comet-type case and let $j$ and $k$ be fixed as
    $k\gg j \geq 0$.
    Suppose that the initial values are chosen as $(X_{0},\dot{X}_{0})
    \in \mathscr{L}_{1}^{(0)}$. 
   There are three undetermined values $\xi_{1},\dot{\xi}_{2},
   \dot{\xi}_{3}$, and three periodicity conditional equations are
   as follows, 
   \begin{align*}
    \tilde{\xi}_{2}(T/4,\xi_{1},\dot{\xi}_{2},
   \dot{\xi}_{3})=0,  \quad 
    \dot{\tilde{\xi}}_{1}(T/4,\xi_{1},\dot{\xi}_{2},
   \dot{\xi}_{3})=0,  \quad
     \dot{\tilde{\xi}}_{3}(T/4,\xi_{1},\dot{\xi}_{2},
   \dot{\xi}_{3})=0 . 
   \end{align*}
 As $T$ is also unknown, the integration time
 can be determined by the orbit's 
 $(k+j+2)$-th passing of the $\xi_{1}\xi_{3}$ plane. 
 The parameter $T$ can be exported as 
 a global variable.
 One can also try to add the period as one undetermined parameter 
 to the conditional equations, like
 $ T/4=T_{0}/4+\xi_{2}(T_{0}/4)/\dot{\xi}_{2}(T_{0}/4)$,
 and $T$ can be calculated iteratively.
 However, it may fail when the precision of $T_{0}$ is poor for the full
 system, or say when the ratio $j/k$ is not small enough.

  A \textsc{SDSP} starts from the $x_{1}x_{3}$
  plane, which is also the $\xi_{1}\xi_{3}$ plane
  at the initial epoch. 
  During the time span $T_{0}/4$,
  as the line of syzygy rotates $(k+0.5)\pi$, 
  the orbit of the infinitesimal body
  hits the $x_{2}x_{3}$ plane perpendicularly
  at the $(j+k+1)$-th passing. Meanwhile,
  the $\xi_{1}\xi_{3}$ plane rotates in the inertial
  frame anticlockwise with the angular velocity 
  $n^{\prime}=1$, and coincides with
  the $x_{2}x_{3}$ plane at $T_{0}/4$.
  The infinitesimal body passes the $\xi_{1}\xi_{3}$ plane
  $j+k+2$ times containing the beginning and the ending
  within one fourth period.
  We use the interpolation method to get the $(k+j+2)$-th
  intesection with the $\xi_{1}\xi_{3}$ plane so as to get
   the periodicity conditional equations. 
   
   Suppose that the initial values are chosen as $(X_{1},\dot{X}_{1})
    \in \mathscr{L}_{2}^{(0)}$. Three parameters $(\tilde{\xi}_{1}$,
    $\tilde{\xi}_{3}$,$\dot{\tilde{\xi}}_{2})$ are needed to be 
    determined. Then the periodicity conditions are as follows, 
   \begin{align*}
    \xi_{2}(T/4,\tilde{\xi}_{1},\tilde{\xi}_{3},
   \dot{\tilde{\xi}}_{2})=0,  \quad 
    \xi_{3}(T/4,\tilde{\xi}_{1},\tilde{\xi}_{3},
   \dot{\tilde{\xi}}_{2})=0,  \quad 
   \dot{\xi}_{1}(T/4,\tilde{\xi}_{1},\tilde{\xi}_{3},
   \dot{\tilde{\xi}}_{2})=0 . 
   \end{align*}
   For example, the $(k+j+2)$ passings can be counted from 
   how many times that $\xi_{2}(t)\xi_{2}(t+\Delta t)\leq 0$ 
   is satisfied, where $\Delta t$ is the integration step.
   And the $(k+j+2)$-th intesection can be precisely acquired if we use
   the Hermite interpolation (or cubic spline interpolation).

   \section{Continuation scheme}
   \label{sec:5}
    Theoretically, the comet- and Hill-type \textsc{SDSP}s are shown
    to exist with the corresponding small parameter arbitrarily small.
   However, there is a lack of numerical evidences
   on calculating these solutions. One may wander to
   know how small the parameter $\varepsilon$ is.
   There is a saying that the small parameter 
   is not necessarily too small in the practical calculations.       
   The continuation scheme is based on Broyden's method with 
   a line-search which can be referred to \citet{Fortran90,XUAS22}.
   Although \citet{Kalantonis03} have pointed out this method,
   we haven't realized their work until very recently.
    Broyden's method with a line-search
  is a global convergent strategy for finding
  roots of the nonlinear equations near the approximate roots. 
  The practical calculations have shown that our algorithm is robust.   
   Let's concisely introduce the necessary formulas 
   and then briefly recall the algorithm.
   
   \subsection{Hermite interpolation}
   Let one solution $\phi(t,X_{0},Y_{0})$ start from
   $\mathscr{L}_{1}^{(0)}$. 
   The periodicity conditions can be derived via the integration
   and the interpolation. We take use of the variable step RKF7(8)
   method to implement the integration and the Hermite interpolation
   to determine the $(k+j+2)$-th intersection with the $\xi_{1}\xi_{3}$
   plane. Suppose we have the $(k+j+2)$-th passing with
   $\xi_{1}\xi_{3}$ plane by the condition 
   $\xi_{2}(t_{1})\xi_{2}(t_{1}+\Delta t)\leq 0$. We want to get the time
   $t=\frac{T}{4}$ which make $\xi_{2}(t)=0$. Linear interpolation
   method has been applied in \citet{Henon03}, however, for a long
   time integration and the double precision, it is better to use a
   smooth interpolation. Here we resort to a vector version of the
   two-point Hermite interpolation for a better 
   precision. Let $t_{2}=t_{1}+\Delta t$, $t\in [t_{1},t_{2}]$,  and let
   \begin{align*}
   &  F_{1}=\xi(t_{1}),  \quad F_{2}=\xi(t_{2}),  \quad
   F_{1}^{\prime}=\dot{\xi}(t_{1}),  \quad 
   F_{2}^{\prime}=\dot{\xi}(t_{2}), 
   \nonumber  \\
   &  l_{1}(t)=\frac{t-t_{2}}{t_{1}-t_{2}}\in[0,1], \quad 
   l_{2}(t)=1-l_{1}(t)=\frac{t-t_{1}}{t_{2}-t_{1}}.
   \end{align*}
   The Hermite interpolation basis functions are
   \begin{align*}
     A_{i}(t)=\left[1-2(t-t_{i})l_{i}^{\prime}(t_{i})\right]l_{i}^{2}(t),
    \quad B_{i}(t)=(t-t_{i})l_{i}^{2}(t), \quad i=1,2.
   \end{align*}
   Our vector form of the Hermite interpolation polynomial is
   \begin{align}
     \textbf{\textit{P}}(l_{1}(t)) = & \sum_{i=1}^{2}\left[F_{i}A_{i}(t)+
     F_{i}^{\prime}B_{i}(t)\right]=
     \left[\Delta t (F_{1}^{\prime}+F_{2}^{\prime})+2(F_{1}-F_{2})
     \right] l_{1}^{3}(t)  \nonumber  \\
     & +\left[\Delta t (-F_{2}^{\prime}-2F_{1}^{\prime})
     +3(F_{2}-F_{1})\right] l_{1}^{2}(t)+ (\Delta t) F_{1}^{\prime}
     l_{1}(t)+F_{1} .
   \end{align}
   The second component of $\textbf{\textit{P}}$ should be equal to zero,
   and the unary cubic equation of $l_{1}$ can be solved analytically.
   The formula of the three complex roots can be referred to \citet{Fortran90} and the real root $l_{1}^{\ast}\in (0,1)$ can be searched. Then the periodicity conditions are approximated by
   $\tilde{\xi}_{2}(l_{1}^{\ast}) \approx 0$, 
   $\dot{\tilde{\xi}}_{1}(l_{1}^{\ast}) \approx 0$, and
   $\dot{\tilde{\xi}}_{3}(l_{1}^{\ast})\approx 0$.
    These periodicity conditions are nonlinear and can be converted
   to a minimum problem. 
   
   \subsection{Broyden's method with a line search}
          
  Consider a general $m$-dimensional
   system of nonlinear equations 
   $\Phi(\textbf{\textit{X}})=0$. The map
   $\Phi: \textbf{\textit{X}}\in \mathbb{R}^{m}
   \rightarrow \mathbb{R}^{m}$ is Lipschitz continuous 
   and differential.
    The objective function is
   \begin{align}
    f_{\texttt{unc}}(\textbf{\textit{X}})=\frac{1}{2} \Phi^{\mathtt{T}}\Phi . 
   \end{align}
   The goal is to linearly search a new point
    \begin{align}
      \textbf{\textit{X}}^{(k+1)}=\textbf{\textit{X}}^{(k)}+\lambda\cdot 
      \textbf{\textit{p}}^{(k)}, \quad \lambda\in (0,1), 
   \end{align}
   starting from the known point $\textbf{\textit{X}}^{(k)}$
   along the descending direction $\textbf{\textit{p}}^{(k)}$,  
   such that $f_{\texttt{unc}}( \textbf{\textit{X}}_{k+1})<
   f_{\texttt{unc}}(\textbf{\textit{X}}_{k})$. The procedure continues until
    $\|\Phi\|_{\infty}$ is smaller than the allowance error $10^{-10}$,
    where the $\infty$-norm represents the maximum component of the
    absolute values. $\textbf{\textit{p}}^{(k)}$ is solved by 
    Broyden's method with the QR decomposition, and $\lambda$
    is solved by the line-search. The fortran subroutine 
    ``broyden(x,check)'' of this method can be found in the book
    \citep{Fortran90}. 
        
   Broyden's method is a kind of quasi-Newton method, and the 
   Jacobian matrix is approximated by $\textbf{\textit{B}}$. By
   Newton's method, we have
   \begin{align}
    \textbf{\textit{B}}^{(k)}\textbf{\textit{p}}^{(k)}=
    - \Phi(\textbf{\textit{X}}^{(k)}).
   \end{align}
   The initial matrix $\textbf{\textit{B}}^{(0)}$ is approximated by
   the forward difference method. If $\textbf{\textit{B}}^{(0)}$ is
   non-singular, $\textbf{\textit{B}}^{(0)}$ can be decomposed
   as $\textbf{\textit{B}}^{(k)}=\textbf{\textit{Q}}^{(k)}
   \textbf{\textit{R}}^{(k)}$ with $k=0$ by the subroutine ``qrdcmp'', 
   where $\textbf{\textit{Q}}$ is 
   the orthogonal matrix and $\textbf{\textit{R}}$ is the upper
   triangular matrix. Then we have the formula
   \begin{align}
    \textbf{\textit{R}}^{(k)}\textbf{\textit{p}}^{(k)}=
    -(\textbf{\textit{Q}}^{(k)})^{\mathtt{T}}
    \Phi(\textbf{\textit{X}}^{(k)}).
   \end{align}
    The matrix $\textbf{\textit{B}}^{(k+1)}$ is updated by
    \begin{align}\label{textbfitBk}
     \textbf{\textit{B}}^{(k+1)}=\textbf{\textit{B}}^{(k)}+
    \textbf{\textit{w}}^{(k)}\otimes \textbf{\textit{s}}^{(k)} ,
    \end{align}
    where 
    \begin{align}
     & \textbf{\textit{w}}^{(k)}=\Phi(\textbf{\textit{X}}^{(k+1)})-
      \Phi(\textbf{\textit{X}}^{(k)})-\lambda \textbf{\textit{B}}^{(k)}
      \textbf{\textit{p}}^{(k)}, \nonumber  \\
     &   \textbf{\textit{s}}^{(k)}= \textbf{\textit{p}}^{(k)}/
     ( \textbf{\textit{p}}^{(k)}\cdot \textbf{\textit{p}}^{(k)}) .
    \end{align}
    The QR decomposition of $\textbf{\textit{B}}^{(k+1)}$ is updated
    by the subroutine ``qrupdt'', which needs $2(m-1)$ Jacobi
    rotations. 
    
     Let's review the backtracking line search strategy. The objective
     function is
     \begin{align}
       \min_{0<\lambda<1}f_{\texttt{unc}}(\textbf{\textit{X}}^{(k)}
       +\lambda\cdot \textbf{\textit{p}}^{(k)})
      =\min_{0<\lambda<1}\varphi(\lambda)  
       \approx \lambda(0)+\varphi^{\prime}(0)\lambda_{1}+
      c_{1}\lambda_{1}^{2} ,
     \end{align}
     where $c_{1}=\varphi(1)-\varphi(0)-\varphi^{\prime}(0)$, 
     and $\lambda_{1}=-\frac{\varphi^{\prime}(0)}{2c_{1}}$.
     In order to calculate the sequent backtracks, the unary cubic equation with two unknowns $c_{2}$ and $c_{3}$, 
     \begin{align}
       \varphi(\lambda)\approx c_{3}\lambda^{3}+c_{2}\lambda^{2}
       +\varphi^{\prime}(0)\lambda+\varphi(0) ,
     \end{align}
     can be solved analytically by giving two different 
     values $\lambda_{2}$ and $\lambda_{3}$. The minimum
     can be solved by $\varphi^{\prime}(\lambda^{\ast})=0$ and let
     $\lambda^{\ast}$ be in $[0.1,0.5]\lambda_{1}$. 
     Redefine $\lambda_{1}=\lambda^{\ast}$ and repeat
     the procedure until $c_{3}$ tends to zero.
     
       There are three key steps for the numerical continuation. 
   Firstly, choose one group of approximate initial values. 
   Secondly, choose the periodicity conditions, which are nonlinear 
   and derived from the numerical integration and interpolation. 
   Thirdly, take use of Broyden's method with a line search to solve
   the nonlinear equations.   
         

     \section{Linear stability}
     \label{sec:6}
     Consider the linear variation of the differential
  equations (\ref{RotDiff}), which can be 
  rewritten as $\dot{\xi}=
  \partial \mathcal{H}^{\texttt{rot}}/\partial
  \eta$, $\dot{\eta}=- \partial \mathcal{H}^{\texttt{rot}}/\partial  \xi$.
  Then the linear variations satisfy
  \begin{align}\label{zeta1Diff1}
   \delta \dot{\xi}=\delta \eta+
   (\delta \xi_{2},-\delta \xi_{1},0)^{\mathtt{T}},
   \quad 
   \delta\dot{\eta}\approx -\frac{\partial^{2} \mathcal{H}
   	^{\texttt{rot}}}{\partial \xi 
   	\partial (\xi,\eta)}
   \left(\begin{array}{c}
   \delta \xi  \\  \delta \eta 
   \end{array}\right), 
  \end{align}
  the formulas can be referred to \citet{Lara02} or one can calculate them by one symbol algebra software.   
 
  Suppose $\left(\xi(T),\eta(T)\right)=\left(\xi(0),\eta(0)\right)$,
  this means that a solution starts from
  $(\xi_{0},\eta_{0})$ and comes back after time $T$. 
  Along this periodic solution, the 
  coefficient matrix $-\frac{\partial^{2} \mathcal{H}
  	^{\texttt{rot}}}{\partial \xi 
  	\partial (\xi,\eta)}$ is periodic. 
  According to the Floquet theorem,
  the stability depends on the eigenvalues of the monodromy matrix,
  which is a fundamental solution matrix
  taking the value after a period. The eigenvalues of 
  $\frac{\partial (\xi(T),\eta(T))}{\partial (\xi(0),\eta(0))}$
  are the same as those of $\frac{\partial (\xi(T),\dot{\xi}(T))}
  {\partial (\xi(0),\dot{\xi}(0))}$ as there exists only elementary matrix 
  transforms between these two monodromy matrices.
    
   Consider $\zeta=
  (\delta \xi_{1},\cdots,\delta \eta_{1},\cdots,
  \delta \eta_{3})^{\mathtt{T}}
  \in \mathbb{R}^{6}$. The differential 
  equation can be derived as (\ref{zeta1Diff1}).   
  The standard fundamental solution matrix about
  $\zeta$ can be denoted as
  $\mathcal{Z}(t)=\mathcal{Z}(t,0)=\left(
    \frac{\partial \zeta_{i}(t,0)}{\partial 
    \zeta_{j}(0)} \right)_{1\leq i,j\leq 6}$
 and satisfies  
  $\mathcal{Z}(0)=\mathbf{I}_{6\times 6}$,  
  where $\mathbf{I}_{6\times 6}$ is a $6\times 6$
  identical matrix, so the monodromy matrix
  $\mathcal{Z}(T)$ can be integrated numerically.
  Some Lemmas can be referred to 
  \citet{MeyerHO}.
  \begin{lemma}
  If $\mathcal{Z}(t)$ is a fundamental solution
  matrix of a linear differential system
  with a $T$-periodic coefficient matrix, then
  $\mathcal{Z}(t+T)=\mathcal{Z}(t)
  \mathcal{Z}(T)$.
  \end{lemma}
  \begin{theorem}
   The fundamental solution matrix $\mathcal{Z}(t)$
   is symplectic for all $t\in \mathbb{R}$.	The
   eigenvalues of $\mathcal{Z}(T)$ (characteristic
   multipliers) are symmetric with respect to
   the real axis and the unit circle.
  \end{theorem}  
  
  \ \ According to the double-symmetry property,
  the $\mathcal{Z}(T)$ can be calculated by
  just integrating $T/4$ \citep{RobinM80}.
  \begin{theorem} \label{theorem3}
  	Let $\mathcal{Z}(t)$ be defined as above.
   	If a \textsc{SDSP} starts from $\mathscr{L}_{1}^{(0)}$,
   	and $\mathcal{Z}(T/4)=\mathcal{Z}_{(0,T/4)}$ is known, then the monodromy matrix of this \textsc{SDSP} 
   	can be calculated by 
   	\begin{align} \label{MonoZT}
   	\mathcal{Z}(T)=\left(\mathcal{N}\left(\frac{T}{2}\right)
   	\mathscr{R}_{1}\right)^{2},
   	\quad \mathcal{N}\left(\frac{T}{2}\right)
   	=\mathcal{Z}\left(\frac{T}{4}\right)
    \mathscr{R}_{2}\mathcal{Z}^{-1}
   \left(\frac{T}{4}\right),
   	\end{align}
   	where $\mathscr{R}_{1}=$ $\mathtt{diag}$ 
   	$(1, -1, -1, -1, 1, 1)$ and $\mathscr{R}_{2}=\mathtt{diag}(
   	1, -1, 1, -1, 1, -1)$. If the solution starts
   	from $\mathscr{L}_{2}^{(0)}$ and 
   	$\tilde{\mathcal{Z}}(T/4)$ is known, then
   	\begin{align*} \label{MonoZT2}
   	\mathcal{Z}(T)=\left(\tilde{\mathcal{N}}(\frac{T}{2}) 	
   	\mathscr{R}_{2}\right)^{2},
   	\quad \tilde{\mathcal{N}}(\frac{T}{2})=
   	\tilde{\mathcal{Z}}(\frac{T}{4})
   	\mathscr{R}_{1}\tilde{\mathcal{Z}}^{-1}(\frac{T}{4}) .
   	\end{align*}
  \end{theorem}
  \begin{proof}
  	According to Liouville's theorem, the 
  	trace of the periodic coefficient matrix
  	equals zero, so the determinant of 
  	$\mathcal{Z}(t)$ equals $1$.
  According to the property of the fundamental
  solution matrix, we have $\mathcal{Z}(T)
  =\mathcal{Z}(T/2)\mathcal{Z}_{(T/2,T)}$.
  For the symmetry $\mathscr{R}_{1}$, we can find that $\mathcal{Z}_{(T/2,T)}
  =\mathscr{R}_{1}\mathcal{Z}^{-1}(T/2)\mathscr{R}_{1}$.
  For the symmetries $\mathscr{R}_{1}$ and $\mathscr{R}_{2}$, we have
  \begin{align*} 
  \mathcal{Z}(T) &= \mathcal{Z}\left(\frac{T}{2}\right)
  \left(\mathscr{R}_{1}
  \mathcal{Z}^{-1}\left(\frac{T}{2}\right)
  \mathscr{R}_{1}\right), \nonumber \\
  \mathcal{Z}\left(\frac{T}{2}\right)
 & =\mathcal{Z}\left(\frac{T}{4}\right)
  \left(\mathscr{R}_{2}
  \mathcal{Z}^{-1}\left(\frac{T}{4}\right)
  \mathscr{R}_{2}\right).
  \end{align*}
  We derive
  \begin{align*} 
  \mathcal{Z}^{-1}\left(\frac{T}{2}\right)
  =\mathscr{R}_{2}\mathcal{Z}\left(\frac{T}{4}\right)
  \mathscr{R}_{2}
  \mathcal{Z}^{-1}\left(\frac{T}{4}\right)
  =\mathscr{R}_{2}\mathcal{N}\left(\frac{T}{2}\right)
  \mathscr{R}_{2}.
  \end{align*}
  Then $\mathcal{Z}(T)$ for the \textsc{SDSP}
  can be calculated by $\mathcal{Z}(T/4)$. The proof
  for the other case of $\mathcal{Z}(T)$
  is similar if the \textsc{SDSP} starts from $\mathscr{L}_{2}^{(0)}$.
  \end{proof}

  The characteristic multipliers of (\ref{zeta1Diff1}) 
  measures the stability of periodic solutions.
  For the three-dimensional \textsc{CRTBP},
  there are six characteristic multipliers which are 
  in pairs of reciprocals. In fact, these
  multipliers can be complex conjugate on a circle
  or can be real reciprocals. 
  Denote one multiplier as $\lambda$, 
  if there exists a $|\lambda|>1$, 
  then the periodic solution is not linearly stable.
  The six multipliers can be represented as
  $\lambda_{i},\frac{1}{\lambda_{i}}$ ($i=1,2,3$). 
  We refer to the well known index for the 
  linear stability in \citet{Henon69,Lara02}, and define
  \begin{align} 
  \rho=\sum_{i=1}^{3}\left(|\lambda_{i}|+
  \frac{1}{|\lambda_{i}|}\right)\geq 6
  \end{align} 
  as the linear stability index for the spacial
  periodic orbits. The trace of $\mathcal{Z}(T)$ can be
  firstly used to judge the stability. If $\rho\geq 
  \mathrm{Tr}(\mathcal{Z}(T))>6$, the orbit is
  linearly unstable. While if $\mathrm{Tr}(\mathcal{Z}(T))\leq 6$, we will
  use $\rho$ to help understand the linear stability. 
  When all the multipliers are on the unit circle,
  one has $\rho=6$, and
  such a periodic solution is linearly stable.  
  In fact, for the numerical 
  errors, we can not get the precise $\rho=6$.
  The orbits with $\rho<6+10^{-4}$ are supposed to 
  be linearly stable in the numerical sense.
  The monodromy matrix $\mathcal{Z}(T)$
  can be calculated by the 
  integration of linear variational equations together with the corresponding differential equations, 
  and by the application of Theorem \ref{theorem3}.
  The eigenvalues are calculated
  by the elemination method and the QR algorithm. 
  The fortran subroutines "elmhes" and "hqr" can be referred to in \citet{Fortran90}. "elmges" 
  transforms a matrix to the upper Hessenberg form,
  and "hqr" solves the eigenvalues by a sequence of 
  Householder transformations.

   \section{Numerical Examples}
   \label{sec:5}
     All the numerical experiments are carried out
  in a linux system
  on a personal computer with 
  Intel Core i7-6500U CPU $@$ 2.50GHz $\times$ 4 
  and $7.7$ GiB memory. If the maximal numerical deviation of the periodicity conditions at one fourth period is no greater than $10^{-9}$, a symmetric periodic orbit is supposed to be found.  
  The time consuming for the continuation can be controled by the maximum number of loops.  
  The symmetric periodic orbits with the designated $k$ and $j$ are found by the optimization method described in section \ref{sec:5}.
  The linear stability of these orbits are estimated
  by the method described in section \ref{sec:6}, and
  the characteristic multipliers are calculated
  by the elemination method and the QR algorithm.
  The fortran subroutines for calculating
  eigenvalues can be referred to 
  "elmhes" and "hqr" in \citet{Fortran90}.

   \subsection{Case $\mu=1/2$}  
   The special case of equal masses in the \textsc{RTBP}
   known as the ``Copenhagen problem". This case
   has been studied with the generalized force potential.
   Here, we confine our experiment on the Newtonian potential.    
   Generally speaking, it is necessary for the perturbation to be small if one wants to continue the 
   periodic orbits from the approximate system
   to the full system.
   In order to estimate the order of the magnitude
   of the perturbation on the approximate system,
   we compare the first-order perturbation term with
   $\mathcal{H}_{0}^{\texttt{rot}}$.
   The first-order perturbation term is
  $\mathcal{H}_{1}^{\texttt{rot}}=
   -(1-\mu)\mu \mathcal{P}_{2}(\xi_{1}/r)/r^{3}$.
    For the nearly circular
   orbits near the infinity, the line of syzygy
   moves much faster than the mean motion
   of the big orbit of the infinitesimal body.
   During one orbit period of the infinitesimal
   body, the averages of $\cos 2\Omega$ and
   $\sin 2\Omega$ equal zeros. We refer to \citet{XbXu19} for the necessary terms and have 
   \begin{align*}
   \frac{1}{2\pi}\int_{0}^{2\pi}\mathcal{P}_{2} 
    (\xi_{1}/r) r^{-3} \mathrm{d}\ell/n 
    =\frac{1}{8}(3\cos^{2}i-1)\Theta^{-3} . 
   \end{align*}   
   It seems that the first-order averaged perturbation equals zero
   when $\cos^{2} i=1/3$ though the first-order perturbation
   does not vanish.   
 
\begin{figure}[h]
  	\center{\includegraphics [scale=1.00]{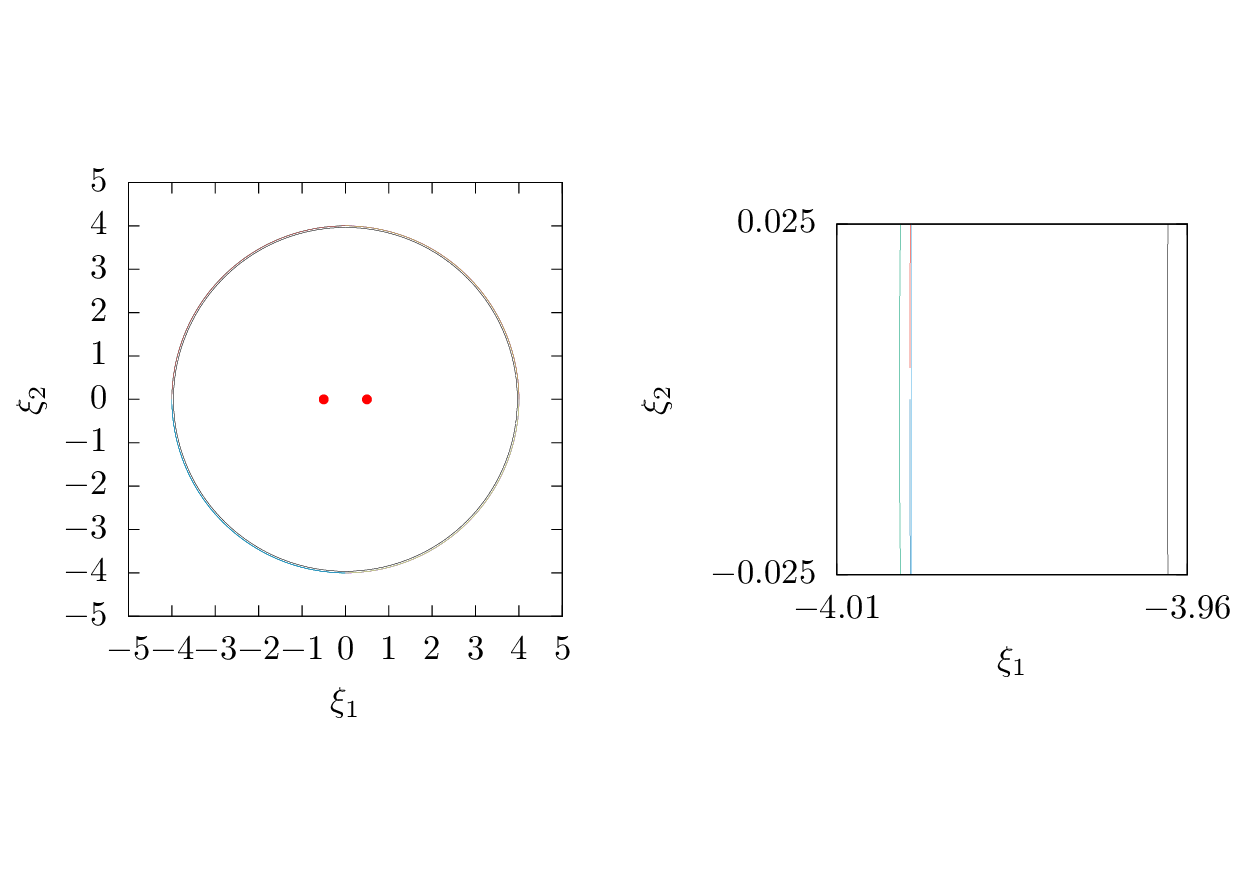}}
  	\caption{Planar symmetric period orbits of family \emph{m} continued from the same initial values in Eq.(\ref{InitialLara}) with $\mu=0.5$. The initial	values of the inner periodic orbit are
  	$(x,y,\dot{x},\dot{y})=(0,3.96199469992294,          
  	4.46677589984367,0)$, and the period is
  	$5.57243120610132$. The outer periodic orbits
  	are continued by the optimization method in this paper. Part of the left-side sketch map is enlarged
  	as the right-side sketch map. Primaries
  are represented by dots. }\label{fig2LaraM}
  \end{figure}

  One set of initial values for the 
  planar \textsc{RTBP} can be referred to
  \citet{Lara02} as follows,
  \begin{align}\label{InitialLara}
 (x_{0},y_{0},z_{0},\dot{x}_{0},\dot{y}_{0},\dot{z}_{0})
 =(0, 4, 0, 4.5, 0, 0)\in\mathscr{L}_{2}^{(0)}, \quad
 T_{0}=5.585 .
  \end{align}
  The periodic orbits continued from this set of values
  belong to family \emph{m}, in which the periodic orbits
  are usually stable.
  Our numerical experiment shows that there are a lot of 
  symmetric periodic orbits near this approximate solution.
  The results achieved by our algorithm are more close
  to the initial values. In Fig.\ref{fig2LaraM},
  there are four near circular retrograde periodic orbits.
  The outer orbit is symmetric with respect to
  the $y$ axis.
  The initial values of the outer orbit is
  $y_{1}=4.00021433614826$, 
  $\dot{x}_{1}=4.50254016243978$, and the period is
  $T_{1}=5.5815257432169$. There are two orbits
  in the middle and almost coincide. 
  One orbit is determined by the double symmetries,
  and the other orbit is determined by 
  the second perpendicular crossing 
  of the positive $y$ axis. The initial values are $(y_{2},\dot{x}_{2})$$=$$(3.99876390419082$, 
  $4.50118225570633)$, and 
  $(y_{3},\dot{x}_{3})$$=$$(3.99857131110253$, 
  $4.50100195228139)$, respectively.
  The periods are 
  $T_{2}=5.5811840743002$ and
  $T_{3}=5.58113868599737$, respectively. 
  The accuracy of the three periodic orbits
  is within $1.1$E-13.

\begin{figure}[h]
  	\center{\includegraphics [scale=1.00]{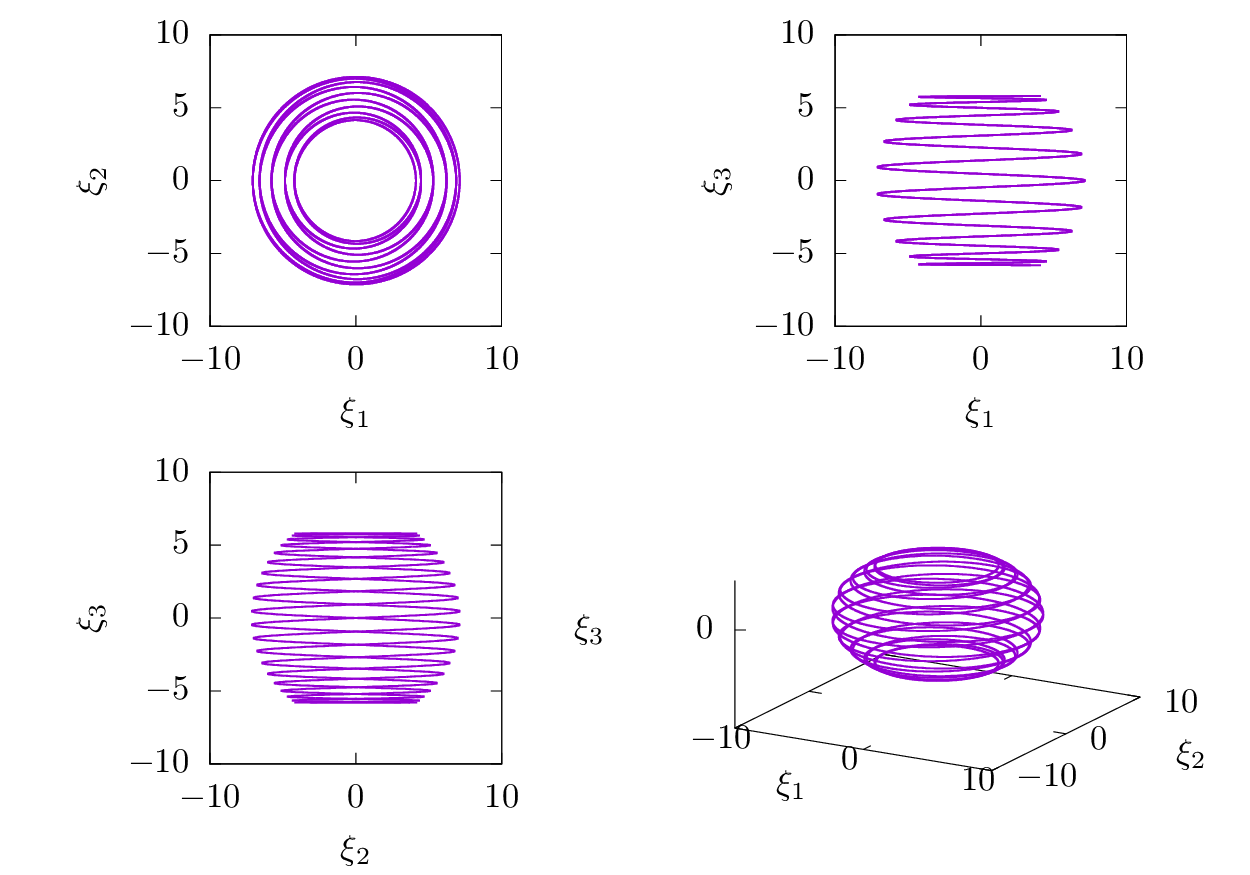}}
  	\caption{An orbit diagram of comet-type \textsc{SDSP}s
  	 in case $(1,+,-,-)$ with
  	$\mu=0.5$, $k=9$, $j=0$, $\cos i=\sqrt{3}/3$.
  	The initial values and the period can be found in Table \ref{Table1}. The characteristic multipliers can be found in Table \ref{Table2}.}\label{fig3Comet}
  \end{figure}

\begin{figure}[h]
  	\center{\includegraphics [scale=0.45]{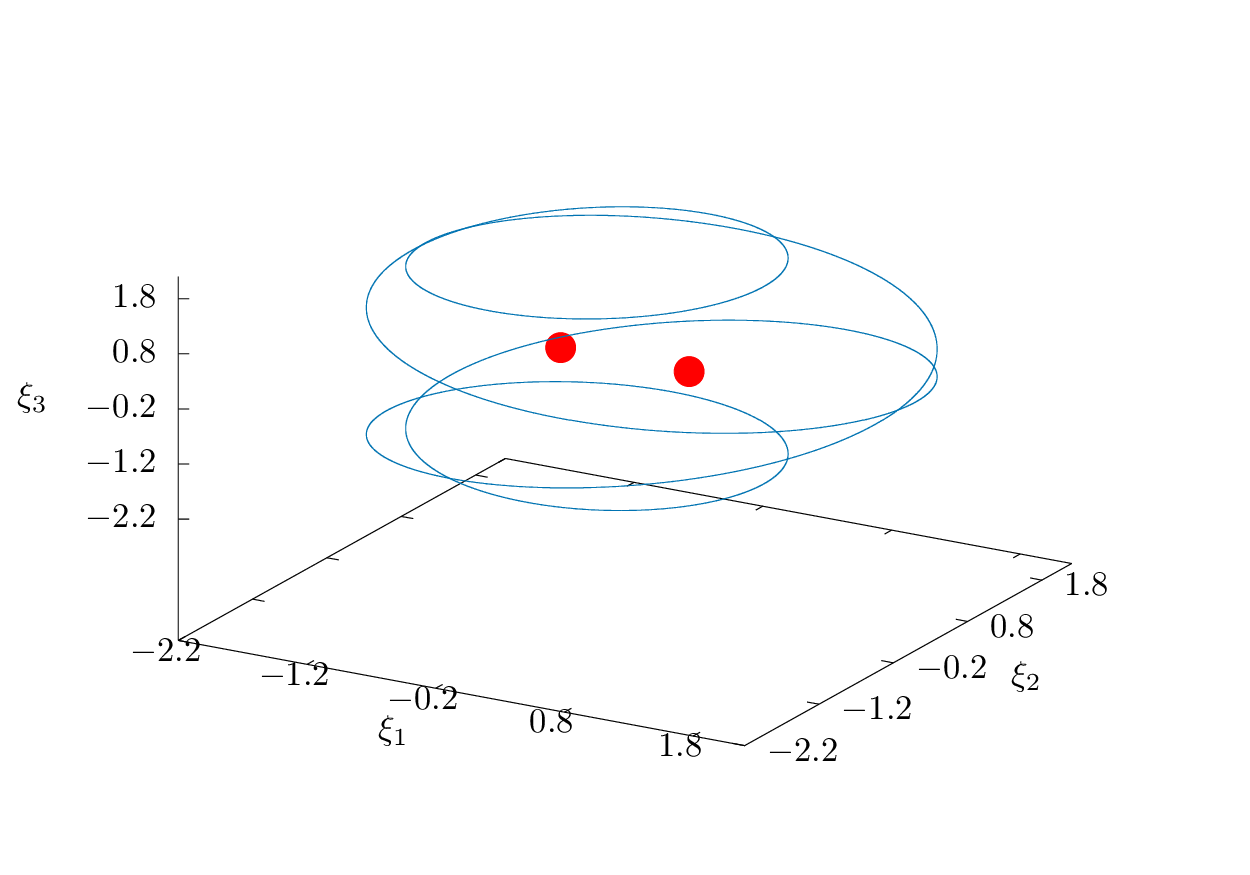}
  		\includegraphics [scale=0.45]{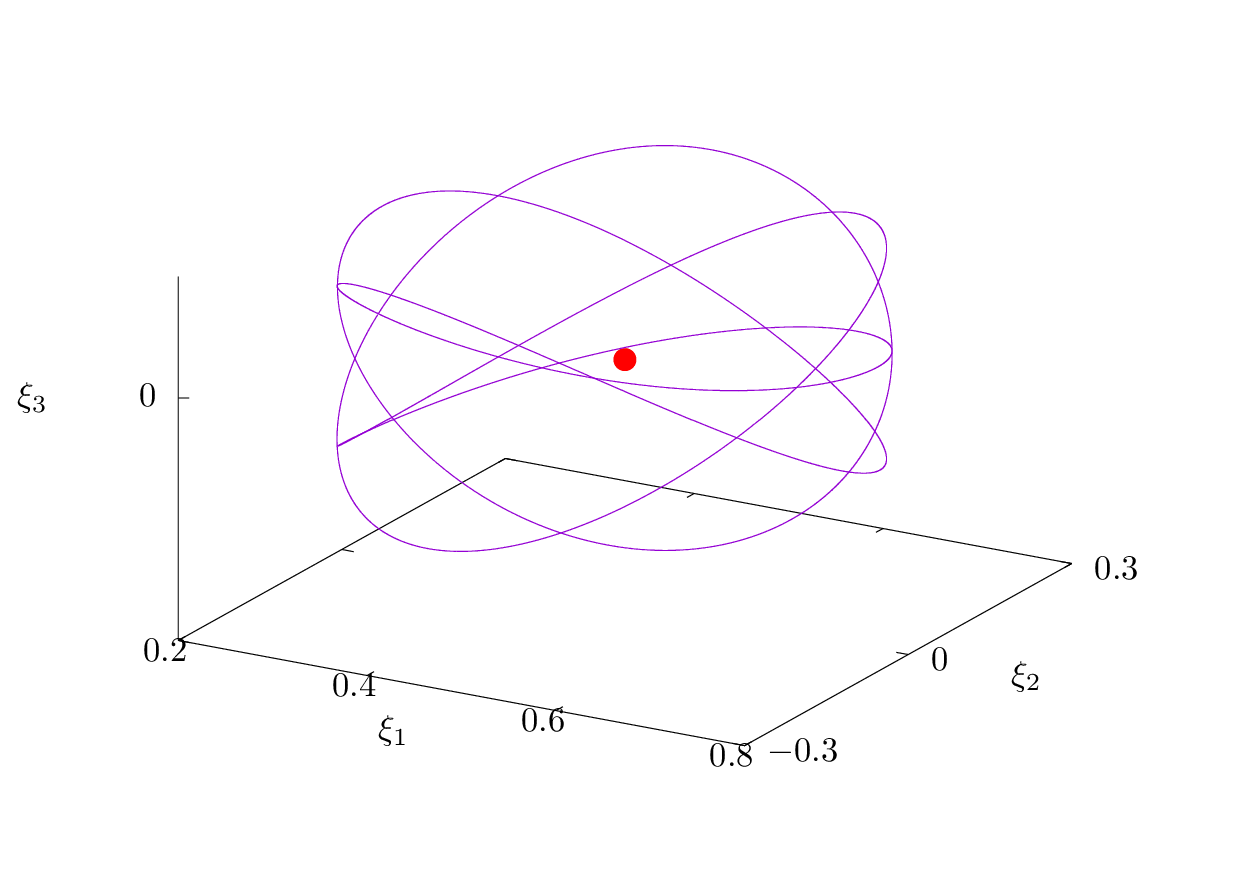}
  		\includegraphics [scale=0.45]{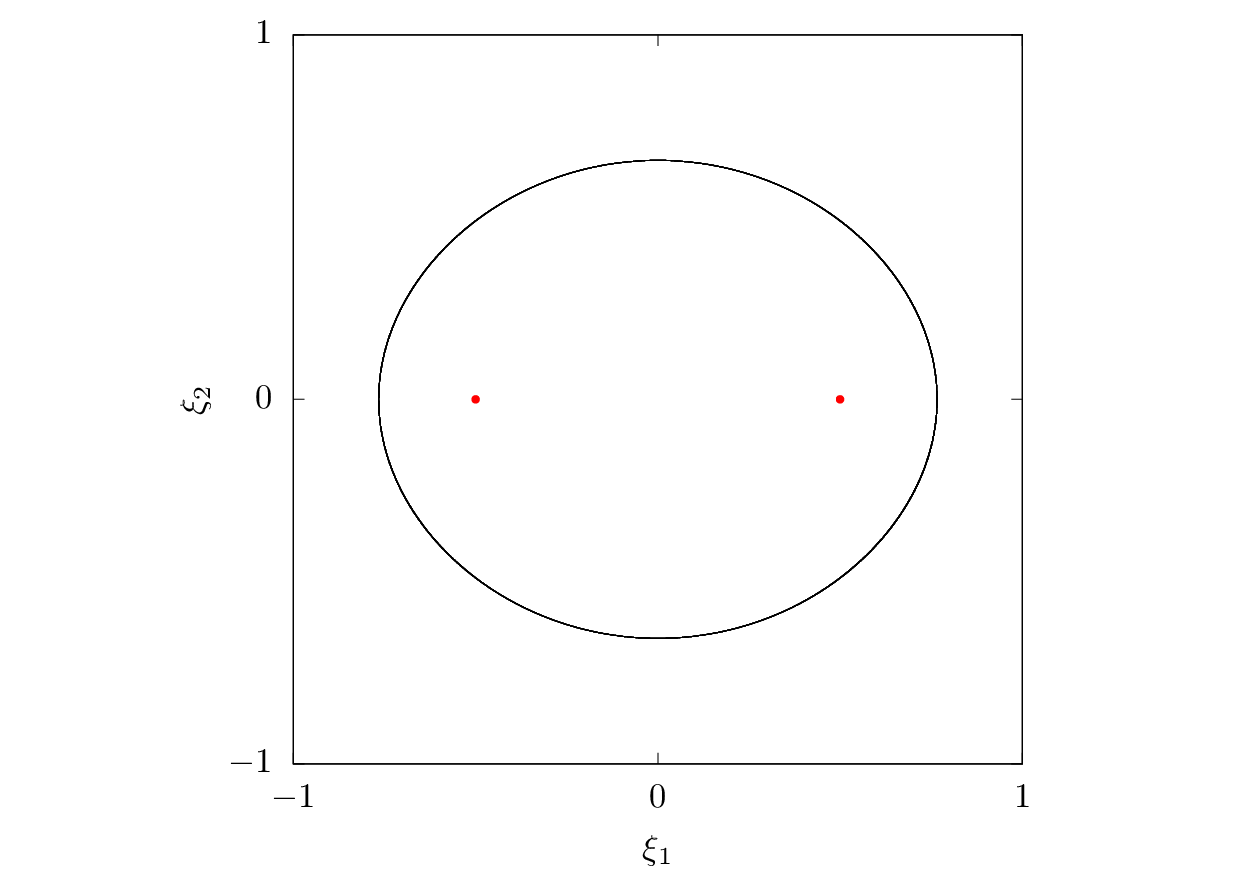}
  		\includegraphics [scale=0.45]{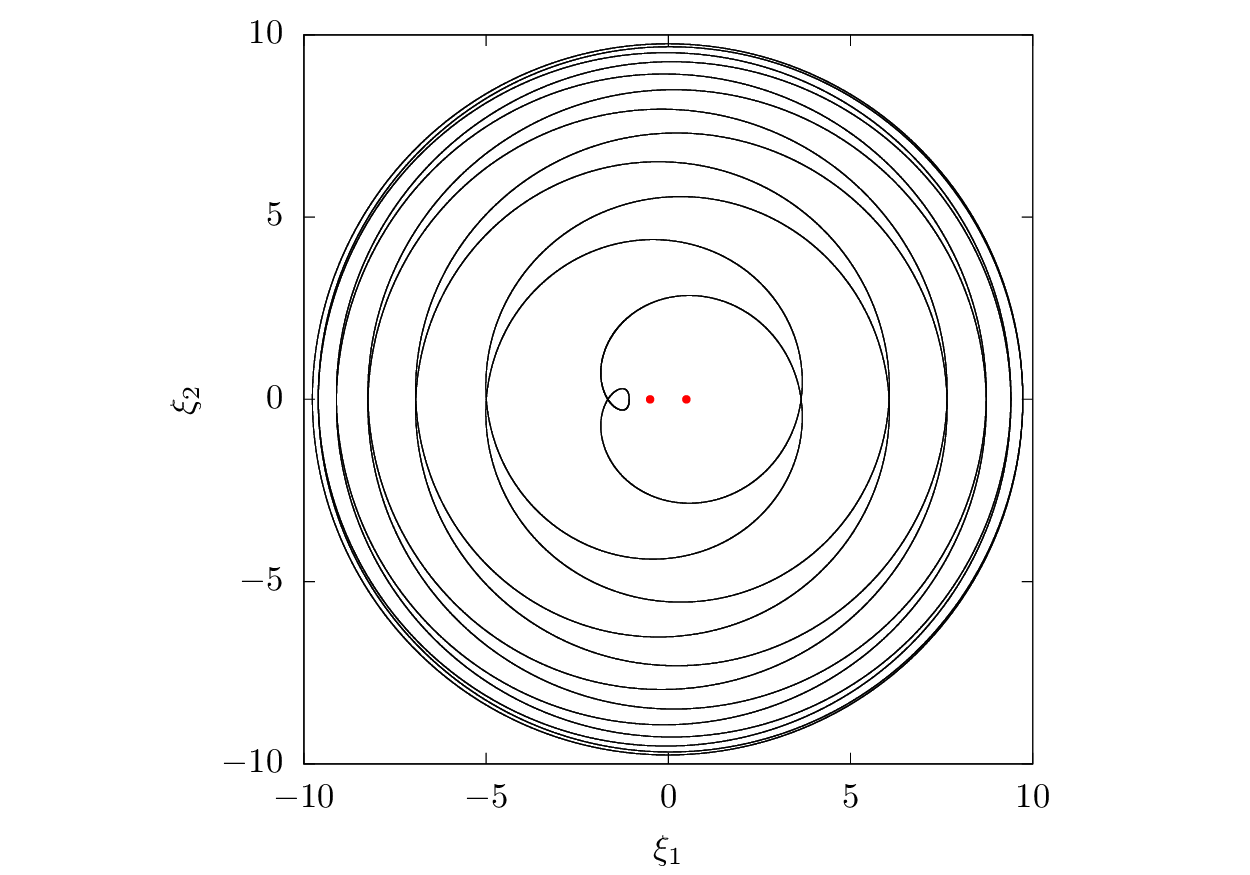}}
  	\caption{Four sketch maps of \textsc{SDSP}s with
  $\mu=0.5$ calculated in Table\ref{Table1}. The upper left one describes the $(1,0)(1,+,-,-)$ type orbit, which is retrograde. The upper right one is continued from the  
  $(1,0)(1,+,-,+)$ type orbit, and the orbit is of 
  Hill-type and retrograde. The lower left one is
  continued from the $(2,0)(1,-,+,-)$ type orbit,
  and the orbit is planar retrograde comet-type.
  The fourth one is interesting and continued from 
  the $(10,0)(1,-,+,-)$ type orbit.
  	}\label{figeqCRTBP4}
  \end{figure}

   For $\mu=0.5$ and $\cos i_{0}=\sqrt{3}/3$, 
  we vary integers $k$, $j$, as well as the directions, 
  and numerically find many \textsc{SDSP}s. 
  Some examples are listed in Table \ref{Table1}.   
  For $k=30,j=0$,
  we get $a_{0}=61^{2/3}\approx 15.496$, $T_{0}=
  30.5\pi\approx 95.819$.
  It takes about several minutes to get the continued
  initial values. For the $(1,+,+,+)$ case, we have   
 $\xi_{1}=15.5061254882711$, $\dot{\xi}_{2}= -15.7370477222493$,    
  $\dot{\xi}_{3}=0.105884508957052$, $T/4=95.81968944276656$.
   The accuracy is within $10^{-10}$ by the program. 
  As we can see, the approximate initial values
  approximate the true values, as the 
  small parameter satisfies
   $\epsilon^{3}=1/61\approx 0.0164$.
   This periodic orbit is nearly
   linear stable, as maximum absolute value
   for the six multipliers is $1.000188$.    
  In order to show robust of the program, a
  bigger $\varepsilon\approx 0.101$ with $k=15,j=0$
  is taken. The corresponding initial values for 
  the \textsc{SDSP}s can be calculated as
 $ a_{0}\approx 9.868$, $n_{0}=1/31$ and
  $i_{0}\approx 0.9553$. More continuation results 
  are listed in Table \ref{Table1},   
  and the corresponding characteristic multipliers are 
  listed in Table \ref{Table2}.
  In order to have a better view of the \textsc{SDSP}s,
  an example for the case $k=9,j=0$ is shown in Figure \ref{fig3Comet}, and four interesting sketch maps
  are shown in Fig.\ref{figeqCRTBP4}.

   
  \begin{table}[!ht]
  	\centering 
  	\tiny
  	\caption{A list of initial values of the 
  	comet-type \textsc{SDSP}s in the $O-q_{1}q_{2}q_{3}$
  	frame with $\mu=0.5$, $\cos^{2} i=1/3 $. }
  	\begin{tabular}{c|cccccc}
  		\toprule
  	$k,j$ & type  & $\xi_{1}$ &  $\dot{\xi}_{2}$  &  $\dot{\xi}_{3}$   &  $T/4$  & 
  	accuracy \\		\hline
  	\multirow{2}{*}{1,0} &  $(1,+,-,-)$ &	2.1188907053948314  &  -2.4745187952972980 &
  	-0.59854164753778971 & 4.7457525451537164 &	3.5E-13 \\
  \multirow{2}{*}{} & $(1,+,-,+)$ & 0.23862606510911777 & 
  	-1.1215624162229199 &  -0.28539427470548040 & 1.4642141631345391 & 3.9E-13 \\  \cline{1-1}
  	\multirow{2}{*}{2,0} & $(1,+,+,+)$ &	3.6836976532989136 & -3.3058283884238149 & 0.36090164760291182 & 10.979823749195759 & 4.2E-13 \\
  	\multirow{2}{*}{} &	$(1,+,-,-)$ & 
  	2.9521280076112233 & -3.2950401216394209 &
  	-0.47588726499054074 &  7.8737465982322021  & 2.8E-12 \\
  	\multirow{2}{*}{} &	$(1,+,-,+)$ & 
  	2.1350003684163883 & -1.6290350406991201 &      -0.45256379929931584 & 13.990486164660886 & 8.9E-13 \\
  	\multirow{2}{*}{} &	$(1,+,+,-)$ & 	
  	2.9509027501751386 & -3.2849171529780969 &       0.48219253123559525 & 7.8732466675287895 & 1.6E-12 \\
  	\multirow{2}{*}{} &	$(1,-,+,-)$ & 
  	0.76592552564005434 & -2.4408679202287282 &      -3.90E-16 & 3.1691491002387817 & 1.9E-12 \\ \cline{1-1}
  \multirow{3}{*}{3,0} & $(1,+,+,+)$ & 3.6772403850959363 & -3.9581936776454136 & 0.44208873366913287 & 11.007264293833867 & 5.3E-12 \\
  \multirow{3}{*}{} & $(1,+,-,-)$ &
  3.6798781095392297  & -3.9832858714520296   &   -0.42694943750882136 & 11.008183184365842 & 6.7E-12  \\
  \multirow{3}{*}{} & $(1,-,-,\pm)$ &
  49.113458003571736 & -49.256155620843664 &        6.057E-15 &  12.529965200235129 & 2.7E-11  \\ \cline{1-1}
  \multirow{2}{*}{4,0} & $(1,+,-,+)$ &
  4.9607122507024686 & -4.6703571931614620 &     -0.34352130501186157 & 17.271018239802125 & 2.3E-12 \\
  \multirow{2}{*}{} & $(1,+,+,-)$ &
  4.3429005143656703 & -4.6167105797228301 &      0.39596630632567387 & 14.146055760667986 & 2.1E-12 \\
  \cline{1-1}
   \multirow{1}{*}{5,0} & $(1,+,+,+)$ &
  	5.5411465631415799 & -5.2750792761853749 &       0.33197989206642536 & 20.414343829320096 & 
  	5.0E-13 \\ \cline{1-1} 	
 \multirow{2}{*}{6,0} & $(1,+,-,-)$ &
 5.5413474283443271 & -5.7885276405967332 &      -0.34650623227584565 & 20.425926276339542 &
 2.6E-12 \\
  \multirow{2}{*}{} & $(1,+,+,-)$ &
  5.5412518156131636 & -5.7874386537688620 &      0.34721382746976359 & 20.425904058280139 & 
  9.2E-13 \\  \cline{1-1}  	
  	\multirow{6}{*}{7,0} & $(1,+,+,+)$ &
  6.6179845229709739  & -6.4253092553046329 &       
  0.33802242435908708 & 26.700202960380690 & 5.1E-12    \\
  	\multirow{6}{*}{} & $(1,+,-,-)$ &
  6.0924064139522818 & -6.3159804821888450 &    
 -0.33860413536729711 & 23.566319867818308 & 2.6E-12 \\
    \multirow{6}{*}{} & $(1,-,+,+)$ & 
 -18.577352612101084 & 18.809427959334485 &        3.261619E-12 &24.822649303391490 & 2.9E-11 \\
  \multirow{6}{*}{} & $(1,-,+,+)$ & -529.21252467445640 &  529.16905512539029 & 9.672974E-12 &
  25.134805803199733 & 5.9E-10  \\
  \multirow{6}{*}{} & $(1,+,-,+)$ & 6.6205322501927721 &
 -6.3942216208109599 & -0.31645096476381246 &
 26.699624375459504 & 8.3E-12 \\
   \multirow{6}{*}{} & $(1,+,+,-)$ &
  6.0935955028860160 & -6.3304575354470680 &
  0.32943323659034535 & 23.566578148878893 & 8.8E-12  \\ \cline{1-1}                             
  	\multirow{2}{*}{8,0} & $(1,+,+,+)$ & 7.1272729885452506 & -6.9295593652543772  &  0.31851635380859106 &
  	29.842069192973479 & 2.7E-12  \\
  	\multirow{2}{*}{}  & $(1,+,-,-)$ & 6.6216917631075791 &   -6.8468408409852222  & -0.31737448415924291 &
  	26.707426285480270 & 3.6E-13  \\     \cline{1-1}  	
  \multirow{2}{*}{9,0} & $(1,+,-,-)$ & 7.1297390619595964  & -7.3463410100233348   
  &  -0.3060359201288402 & 29.848479121555926 & 8.7E-12 \\  
  \multirow{2}{*}{} & $(1,+,+,-)$ & 7.1297376094757015 &      -7.3463221723251877 &      0.30604823680097659 &
  29.848478854086675 & 1.7E-12  \\  \cline{1-1}
  \multirow{2}{*}{10,0} & $(1,+,-,-)$ & 
   8.5583850805287351 & -8.7666166718563456 & 
 -0.27142494617193946 & 39.272356628883159 & 2.2E-12 \\
    \multirow{2}{*}{} & $(1,-,+,-)$ &  
  -1.0730217078285889 & -0.49407005211052968 &  -5.45E-16 & 41.614859822995612 & 5.6E-12 \\  \cline{1-1}  
  11,0 & $(1,+,+,+)$ & 8.5592521023082604 &  -8.3314029327793140 & 0.25512483678988157 & 39.26722889864530 & 2.1E-13 \\ \cline{1-1}
   12,0 & $(1,+,+,-)$ & 8.5575443252626222  &     -8.7549518653594625 & 0.27940197124395649 &
   39.272229816054946  & 2.6E-12 \\ \cline{1-1}   
  \multirow{3}{*}{13,0} & $(1,+,-,-)$ &
   9.0072748204920021 & -9.1999533787013306 &     -0.27212029515012948 & 42.413599374304304 & 1.2E-11 \\
   \multirow{3}{*}{} & $(1,+,+,-)$ & 9.0072776203970779 & 
 -9.1999941511817198 & 0.27209338579659731 &
  42.413599786538946 & 9.1E-12 \\
   \multirow{3}{*}{} & $(1,-,+,-)$ &
  14.438813086606565 & -14.702104461048878 &       -5.9481E-14 & 43.194656094698310 &
  2.2E-11 \\ \cline{1-1}
  \multirow{4}{*}{14,0} & $(1,+,-,-)$ &
  9.4460814515453766 & -9.6347244992809653 &     -0.26534929964099480 & 45.555006655227537 &
  4.5E-12 \\
  \multirow{4}{*}{} & $(1,-,-,\pm)$ &  0.28850044822865212 & -1.3354492959420852 & -5.00904419258058E-3 &
  7.4555703223513472 & 3.0E-12 \\  
   \multirow{4}{*}{} & $(1,+,-,+)$ & 9.8736749812586808 & -9.7007288800802325 & -0.26732964483165583 &
   48.693044346780553 & 3.3E-12 \\
   \multirow{4}{*}{} & $(1,+,+,-)$ &
 9.4460700748764115 & -9.6345552350606347 &      0.26546154518076970 & 45.555005058614071 &
 9.5E-12 \\     \cline{1-1}
 \multirow{4}{*}{15,0} & $(1,+,+,+)$ & 10.2932688959135  
 & -10.12713919866396 &  0.2638860082896899 & 51.834797140512 & 2.6E-11 \\    
 \multirow{4}{*}{} & $(1,+,-,-)$ & 9.87485373076774 
  & -10.0586486565369 &   -0.260003593241237 
  & 48.6964301233740  & 1.7E-11  \\   
  	\multirow{4}{*}{} &
  $(1,+,+,-)$ & 9.87483813928410 & -10.0584106945163  
  & 0.260160733674865   &  48.6964280189064  
  & 2.5E-11   \\   
  	\multirow{4}{*}{} &
  $(1,-,+,-)$ & 111.166965114126 & -111.072119901553 & 
  -1.65E-16 & 50.308404483550	& 2.3E-10  \\  \cline{1-1}
  \multirow{2}{*}{16,0} & $(1,+,-,-)$ &
  10.294565680577234 & -10.474446524271785 &    -0.25472084351795954 & 51.837882368188239 &
  1.5E-11 \\  
   \multirow{2}{*}{} & $(1,-,+,-)$ &  
 -130.69165443604271 & 130.77912832998516 &
 -3.8E-21 & 53.371352860464135 & 8.7E-10 \\ \cline{1-1}
 \multirow{2}{*}{2,1} & $(1,+,+,-)$ &
 1.5398777196321236  & -2.1003537437909281 &      0.60576718932978935 & 8.1243671768449133 & 1.5E-12 \\
  \multirow{2}{*}{} & $(1,-,+,-)$ &
  -12.403494408752263 & 12.119388455806499  &     -8.48027E-13 & 12.860963770383160 & 8.1E-13\\ \cline{1-1}
  \multirow{2}{*}{3,1} & $(1,+,+,+)$ &
  1.6885402394246654 & -1.2610261655074169 &      0.61915290612874152 & 26.345073824037087 & 
  6.3E-12 \\
  \multirow{2}{*}{} & $(1,+,-,-)$ &
  1.8123445656089838 & -2.2376204712686123 &     -0.62768878902348968 & 11.145122123904219 & 5.3E-12
  \\  \cline{1-1}
  \multirow{1}{*}{4,1} & $(1,+,-,+)$ & 2.2892810343583170  
  & -1.9105060688172568 & -0.53973858377315242 &
  26.538813770714953 & 2.4E-12 \\ \cline{1-1}
  \multirow{2}{*}{27,1} & $(1,-,+,+)$ &
  7.9391814708572346 & -8.2143772561892163 &
  -7.8E-19  & 84.768944749702342  & 2.8E-11 \\
   \multirow{2}{*}{ } & $(1,+,+,-)$ &
   6.9624084467162106 & -7.1812331756665850 &      0.30997572793971845 & 86.404315493643836 & 1.4E-11
   \\ \cline{1-1}
  \multirow{2}{*}{28,1} & $(1,+,-,-)$ &
 7.1297594769831187 & -7.3466056451867638 &     -0.30586274368110572 & 89.545448635525830 & 6.0E-11 \\  
   \multirow{2}{*}{} & $(1,+,+,-)$ & 7.1297457413286169 & -7.3464276225729277 & 0.30597927180014994  &
  89.545441053569306 & 1.4E-10 \\   
  		\bottomrule
  	\end{tabular}\label{Table1}
  \end{table}

 \begin{table}[!ht]
     	\centering 
     	\tiny
 \caption{Characteristic multipliers of the 
 comet-type \textsc{SDSP}s in the Table \ref{Table1}. }
     	\begin{tabular}{c|cccc}
     	\toprule
  $k,j$ & type  & $\lambda_{1}$ &  $\lambda_{2}$  &  $\lambda_{3}$     \\		\hline
  \multirow{2}{*}{1,0} &  $(1,+,-,-)$ &	1.102364
  -i 6.6E-7  &  1.011916 & 1.000000   \\ 
  \multirow{2}{*}{} & $(1,+,-,+)$ & -0.4326865-i 0.901544 & 
    1.016881 &  1.000000  \\  \cline{1-1} 
   \multirow{2}{*}{2,0} & $(1,+,+,+)$ &	
   0.997589-i 6.9404E-2 & 1.000007 & 1.000610 \\ 
    \multirow{2}{*}{} &	$(1,+,-,-)$ & 
  0.998991-i 4.4919E-2   &	1.000002 & 1.000610 \\ 
    \multirow{2}{*}{} &	$(1,+,-,+)$ & 
     1.746796 & 1.000001 &    1.353232 \\  
    \multirow{2}{*}{} &	$(1,+,+,-)$ & 	
     0.999239-i 3.9016E-2 & 1.000001 & 1.000667 \\ 
    \multirow{2}{*}{} &	$(1,-,+,-)$ & 1.000117 & -0.692136-i 0.721767 & 57016.27 \\  \cline{1-1} 
    \multirow{3}{*}{3,0} & $(1,+,+,+)$ & 0.999835-i 1.8155E-2 & 1.000000-i 1.86E-5 & 1.000002 \\ 
    \multirow{3}{*}{} & $(1,+,-,-)$ & 0.999593-i 2.8541E-2  & 1.000000-i 4.0E-7  & 1.000027 \\ 
    \multirow{3}{*}{} & $(1,-,-,\pm)$ &
   0.989418-i 0.145096  & 1.000000-i 4.7E-6 &        0.989414-i 0.145119 \\ \cline{1-1} 
    \multirow{2}{*}{4,0} & $(1,+,-,+)$ &
     0.999666-i 2.5831E-2 & 1.000000-i 1.26E-5 &     1.000000-i 3.3E-7  
     \\ 
    \multirow{2}{*}{} & $(1,+,+,-)$ &
 0.999822-i 1.8875E-2 & 1.000000-i 1.47E-5 & 1.000001  
  \\  \cline{1-1}  
    \multirow{1}{*}{5,0} & $(1,+,+,+)$ &
 0.999833-i 1.8269E-2 & 1.000000-i 1.20E-5 & 1.000000      \\ \cline{1-1} 	
 \multirow{2}{*}{6,0} & $(1,+,-,-)$ & 0.999915-i 1.3040E-2 & 1.000017 & 1.000000-i 1.0E-6 \\ 
    \multirow{2}{*}{} & $(1,+,+,-)$ &
  0.999918-i 1.2780E-2 & 1.000000-i 1.23E-5 &    
  1.000000-i 3.3E-7  \\  \cline{1-1} 	
  \multirow{6}{*}{7,0} & $(1,+,+,+)$ &
  0.999995-i 3.0117E-3  & 1.000000-i 1.33E-5 &       
     		1.000001  \\ 
  \multirow{6}{*}{} & $(1,+,-,-)$ & 
  0.999969-i 8.1407E-3 & 1.000024 &    
     		1.000000-i 1.9E-6 \\ 
    \multirow{6}{*}{} & $(1,-,+,+)$ & 
    1.000000-i 1.47E-5 & 0.3250865-i 0.945684 & 
    0.323810-i 0.946122  \\  
    \multirow{6}{*}{} & $(1,-,+,+)$ & 
    0.999966-i 8.2582E-3 &  1.000000-i 9.6E-6 & 
    0.999966-i 8.2582E-3 \\ 
    \multirow{6}{*}{} & $(1,+,-,+)$ & 0.999957-i 9.2678E-3 & 1.000025 & 1.000000-i 1.6E-6 \\ 
    \multirow{6}{*}{} & $(1,+,+,-)$ &
    0.999939-i 1.1087E-2 & 1.000000-i 4.79E-5 &
     1.000000 \\ \cline{1-1}  	                     
    \multirow{2}{*}{8,0} & $(1,+,+,+)$ & 
    0.999990-i 4.5077E-3 & 1.000012  & 
     1.000001  \\ 
    \multirow{2}{*}{}  & $(1,+,-,-)$ & 
    0.999960-i 8.9988E-3 &  1.000019  & 
    1.000000-i 3.9E-7 \\ \cline{1-1}  
    \multirow{2}{*}{9,0} & $(1,+,-,-)$ & 
  0.999970-i 7.7102E-3 & 1.000033 
  &  1.000000-i 1.7E-6 \\    
  \multirow{2}{*}{} & $(1,+,+,-)$ & 
  0.999970-i 7.7071E-3 & 1.000022 & 
  1.000000-i 4.6E-7 \\  \cline{1-1}  
  \multirow{2}{*}{10,0} & $(1,+,-,-)$ & 
  0.999977-i 6.8389E-3 & 1.000000-i 3.82E-5 & 
     		1.000000-i 5.1E-7 \\  
  \multirow{2}{*}{} & $(1,-,+,-)$ &  
   34130636.7 & 1.000000 &  
   1.1072E-2-i 0.999939 \\  \cline{1-1}   
   11,0 & $(1,+,+,+)$ & 0.999952-i 9.7881E-3 &  
  1.000000+i 4.08E-5 & 1.000000 \\ \cline{1-1}  
  12,0 & $(1,+,+,-)$ & 0.999986-i 5.3307E-3  &
  1.000045   &  1.000000   \\ \cline{1-1}   
  \multirow{3}{*}{13,0} & $(1,+,-,-)$ &
   0.999988-i 4.8483E-3 & 1.000040 & 1.000000-i 1.4E-6  \\
  \multirow{3}{*}{} & $(1,+,+,-)$ & 0.999988-i 4.8531E-3 & 
1.000023 & 1.000000-i 4.1E-7  \\
   \multirow{3}{*}{} & $(1,-,+,-)$ & 1.000050 & -0.99998-i 6.1329E-3 &  -0.999930-i 1.1807E-2 \\ \cline{1-1}
 \multirow{4}{*}{14,0} & $(1,+,-,-)$ &
0.999990-i 4.4689E-3 & 1.000012 & 1.000000-i 6.0E-7 \\
 \multirow{4}{*}{} & $(1,-,-,\pm)$ &
 0.642881-i 0.765966 & 1.000108 & 1.000000 \\  
 \multirow{4}{*}{} & $(1,+,-,+)$ & 0.999996-i 2.8649E-3 & 
 1.000023 & 1.000000-i 6.3E-7 \\
  \multirow{4}{*}{} & $(1,+,+,-)$ &
  0.999990-i 4.4504E-3 & 1.000000-i 2.1E-7 & 1.000046  \\     \cline{1-1}
 \multirow{4}{*}{15,0} & $(1,+,+,+)$ &  0.999997-i 2.3250E-3  & 1.000000-i 4.26E-5 & 1.000000 \\    
 \multirow{4}{*}{} & $(1,+,-,-)$ & 0.999992-i 4.0138E-3 
  & 1.000000-i 4.39E-5 & 1.000000 \\   
 \multirow{4}{*}{} &
 $(1,+,+,-)$ & 0.999992-i 3.9895E-3 & 1.000000-i 5.70E-5  
  & 1.000001   \\   
 \multirow{4}{*}{} & $(1,-,+,-)$ & 0.985298-i 0.170843 & 1.000000-i 5.08E-5 & 0.985297-i 0.170848  \\  \cline{1-1}
 \multirow{2}{*}{16,0} & $(1,+,-,-)$ &
 0.999993-i 3.6817E-3 & 1.000000-i 4.34E-5 & 1.000001 \\  
 \multirow{2}{*}{} & $(1,-,+,-)$ &  
  0.989809-i 0.142402 & 1.000025 & 0.989808-i 0.142405  \\ \cline{1-1}
  \multirow{2}{*}{2,1}  & $(1,+,+,-)$ &
   0.393416-i 0.919360  & 1.062582 & 1.000000-i 3.8E-6 \\
  \multirow{2}{*}{} & $(1,-,+,-)$ &
  0.383759-i 0.923433 & 1.000000+i 9.6E-6   &     0.381100-i 0.924534  \\ \cline{1-1}
  \multirow{2}{*}{3,1} & $(1,+,+,+)$ &
   2.97908E+9 & 17.5693 &  1.000651 \\
  \multirow{2}{*}{} & $(1,+,-,-)$ &
  0.993595-i 0.112999 & 1.008120 & 1.000000-i 6.1E-6 
  \\  \cline{1-1}
  \multirow{1}{*}{4,1} & $(1,+,-,+)$ & 
  0.989916-i 0.141655 & 1.000019 & 1.000000-i 5.3E-6  \\ \cline{1-1}
  \multirow{2}{*}{27,1} & $(1,-,+,+)$ &
  0.999999-i 1.0386E-3 & 1.000012 & 0.908355-i 0.418201  \\
  \multirow{2}{*}{ } & $(1,+,+,-)$ &
   0.999711-i 2.4037E-2 & 1.000000-i 2.84E-5 &  1.000005 \\ \cline{1-1}
  \multirow{2}{*}{28,1} & $(1,+,-,-)$ &
 0.999732-i 2.3172E-2 & 1.000000-i 6.78E-5 & 
 1.000000-i 2.1E-6  \\  
   \multirow{2}{*}{} & $(1,+,+,-)$ & 
   0.999729-i 2.3260E-2 & 1.000034 & 1.000000-i 2.7E-6 \\   
     		\bottomrule
     	\end{tabular}\label{Table2}
     \end{table}

  \subsection{Sun-Jupiter case}
\begin{figure}[h]
     	\center{\includegraphics    		[scale=1.00]{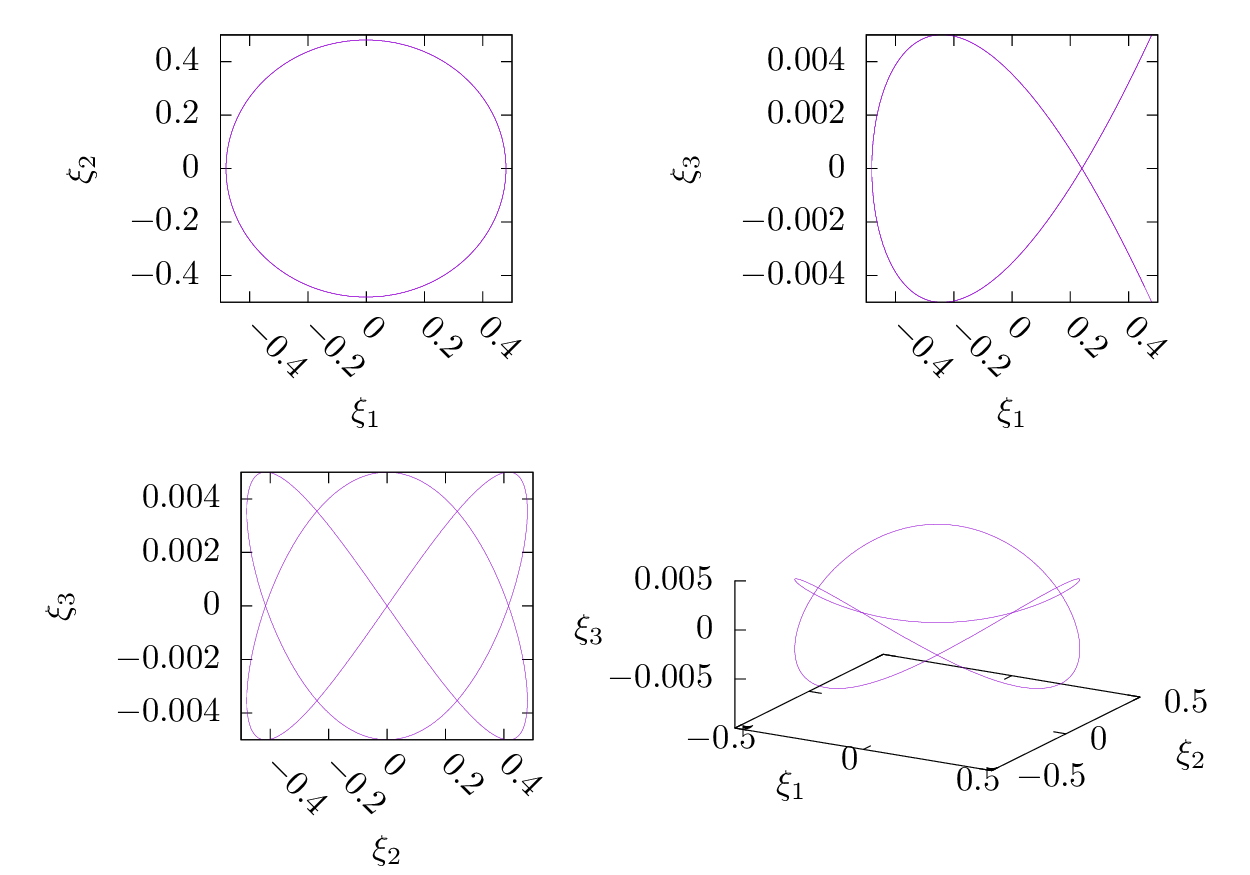}}
     	\caption{An unstable \textsc{SDSP} in the 
  \emph{i}$_{3v}$ family found by \citet{Kazantis79}.
  The orbital information includes $1-\mu=0.00095$,
     		$\tilde{\xi}_{01}\approx$0.47926860,   $\tilde{\xi}_{03}\approx$0.005,
 $\dot{\tilde{\xi}}_{2}$$\approx$0.96274159, 
 $T/4\approx$1.56983, and the Hamiltonian $\mathcal{H}^{\texttt{rot}}\approx -1.733535$.   }\label{figIv3Ex1}
     \end{figure}

  The mass of the Sun is about 1.988500E+30 \texttt{kg},
  while the mass of the Jupiter is about
  1.89813E+27 \texttt{kg}, so we derive the mass
  ratio $1-\mu\approx$ 9.5364E-4.
  However, if the masses of the satellites and the rings
  are included in the Jupiter system, the mass ratio
  $1-\mu$ may approximate 9.5388E-4 according to 
  \citet{BrunoV06}.    
  Consider the calculation of the Hill-type \textsc{SDSP}s
  around the Sun. 
  We use our numerical scheme to check some
  \textsc{SDSP}s calculated by \citet{Kazantis79}. 
  One example is shown in Fig.\ref{figIv3Ex1}.
  The orbit is found to be in the $O-q_{1}q_{2}q_{3}$
  frame, and the accuracy of the initial values given in 
  \citet{Kazantis79} is of order $10^{-4}$ at a fourth
  period. This orbit has
  two nearly perpendicular crossing to the $q_{2}$ axis
  during the time $T/4$. Our numerical experiment with
  $k=j=0$ convinced this result by the accuracy
  $10^{-6}$. However, this \textsc{SDSP} is 
  complex unstable, as
  the index $\rho\approx$ 2.0E+18.
  Our improved initial information is    
     $\tilde{\xi}_{1}=$0.47926856385,
     $\tilde{\xi}_{3}=$5.000000496E-3,
     $\dot{\tilde{\xi}}_{2}=$0.962741651771, and
     $T=$1.569831796427.
 In addition, some more examples are supplied
 in Table \ref{Table3}.   
 
\begin{figure}[h]
     	\center{\includegraphics [scale=1.00]{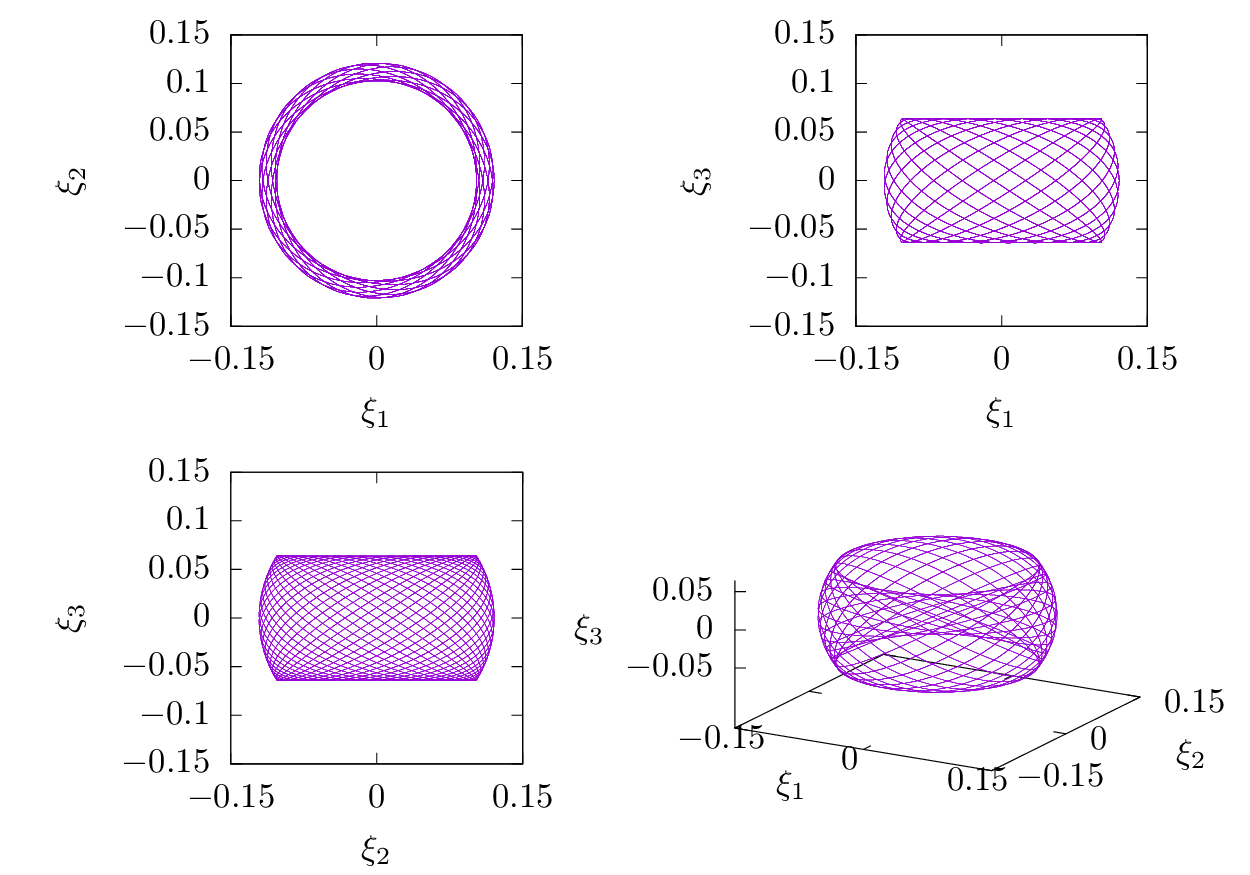}}
     	\caption{One example of Hill-type \textsc{SDSP}
     		around $m_{2}$ in the case $(1,+,+,+)$
     		with $m_{2}=\mu=0.06$, $k=0,j=10$, $\cos i=\sqrt{2}/2$ is supplied. 
 The integration time span is $20$ and its period $T$ is 
 $4$ times of $1.568164286834137$. The initial values 
$(\xi_{1}$,$\dot{\xi}_{2}$,$\dot{\xi}_{3})$ are
$0.12098046779638547$, $2.24272842162778$E-3,$1.4780876804155$E-3, respectively.
     	}\label{figHill006}
     \end{figure}

   If we fix $\mu=0.06$ and $k=0,j=10$,
 then the small parameter $\varepsilon$ satisfies $\varepsilon^{3}=1/21$.      
 Let $a_{1}=\varepsilon^{2}=
 1/(21)^{2/3}\approx 0.13138$, then $|n|=\mu^{1/2} \varepsilon^{-3}
 =21 \sqrt{\mu}$. There is no restriction on the inclination $i$, 
 so we can suppose $i_{1}=\pi/4$. The \textsc{SDSP}
 under this case can be calculated 
 in Fig.\ref{figHill006}. The orbit is of
 multible revolutions and linearly stable.
 The Hill-type \textsc{SDSP}s in Hill's lunar
 problem can also calculated, as many examples
 are shown in Table \ref{Table4}.


    \begin{table}
   	\centering 
   	\tiny
   	\caption{ A list of initial information of the
   	\textsc{SDSP}s in the $m_{2}-q_{1}q_{2}q_{3}$ frame with 
   	$m_{1}=1-\mu=0.00095388$ and $\cos i=\sqrt{2}/2$.}
   	\begin{tabular}{c|cccccl}
   		\toprule
 (k,j) type  & $q_{1}$ &  $\dot{q}_{2}$  &  $\dot{q}_{3}$   &  $T/4$  & accuracy & $\rho$  \\
   		\hline
 (0,1) $(1,+,+,+)$ & 0.34184419200192950 &       0.57007838000595457 &  1.4462000467551235 &
 1.5706863145480114 & 4.7E-13 & 6.00057  \\
 $(1,+,-,-)$ & 0.48068891829543647 &      -1.4998370691818192  & -1.0198406921902150 &
 1.5710001462672483 & 1.5E-13 & 6.00329  \\ 
 $(1,+,-,+)$ &	0.48068899192513403  &     -1.5009002236403080  &  1.0187770988239691 &    
 1.5710006641396190 & 2.2E-13 & 6.00328 \\
 $(1,-,-,+)$ & -0.20851692493061830 &      -2.7837590331097442 & -3.2E-19 & 6.2822221595431698 &
 3.9E-13 & 6.59924  \\
  (0,2) $(1,+,-,-)$ &  0.34192369320122085  &     -1.5508057561554771  &  -1.2085490060598960 &
  1.5709567465854295 & 1.9E-13 & 6.00008 \\
  $(1,-,+,+)$ & -0.34189556538145993 &       1.5412261758084023 &  1.2181591034147481 &     
  1.5709581527177769 & 1.2E-13 & 6.00000 \\
  $(1,+,-,+)$ & 0.34192369381559801 &      -1.5506518035876165 & 1.2087029678246444 &
  1.5709567175162318 & 4.5E-12 & 6.00008  \\       
 $(1,+,+,-)$ &  0.48036437204322080  &     0.58745648907752490   &   -0.97014033755774143  &   
 4.7118144439134388  & 1.2E-12 & 6.12459  \\
  $(1,-,+,-)$ &  -0.34187962005576622 &      0.94062898595117284 &  -1.6012449725489522 &
  1.5708690424098430 &  2.2E-12 & 6.00000 \\  
 (0,3) $(1,+,+,+)$ & 0.23103375232088796 &       1.2532539777939462  &  1.4564805002024286 &
  1.5707061065760803 & 5.4E-13 &  6.00000 \\
  $(1,+,-,-)$ & 0.27320731430162748 &      -1.6256446985588056  &  -1.3519368577710080 &
  1.5709114612538715 & 1.2E-12 & 6.00000  \\
  $(1,+,+,-)$ & 0.23103364333916993 &       1.4004683532340454 & -1.2894297512200859 &
  1.5706952744331948 & 3.3E-12 & 6.00000 \\     
   $(1,-,+,-)$ & -0.27319062603190147 &        1.3338747568633496  & -1.5912567702349809  &   
   1.5708826068536563 & 4.0E-13 & 6.00000 \\   
 (0,4) $(1,+,-,-)$ & 0.23105812269605033  &     -1.7014380953293331 & -1.4703157394679107 &  1.5708855791311456  & 1.2E-13 & 6.00000  \\
 $(1,+,+,-)$ & 0.20210680944323420 &       1.4506267251270490 &  -1.4872032053190620 &     
 1.5707191990654461 & 4.7E-12 & 6.00000 \\ 
 (0,5) $(1,-,-,-)$ & -0.18080434996855890 &      -2.1341178956913809 &  -0.40844779246719487 &     
 1.5707060448070356 & 4.4E-13 & 6.00000 \\
 $(1,+,+,-)$ & 0.18080794122631860 &       1.5413471550336371 & -1.5999217329112894 &
 1.5707323042090662 & 6.8E-13 & 6.00000 \\
 $(1,+,-,+)$ & 0.20212338455189643 &       -1.7747477348967273 & 1.5715116024107758 &
 1.5708692178803645 & 4.4E-12 &  6.00000 \\
 $(1,+,-,-)$ & 11.073693138687092 &  -11.373259558709655 
 & -2.87518E-16 & 18.351486458331845 & 6.1E-12 & 6.00000 \\
   $(1,-,+,-)$ & -0.20211773027265062 &       1.4577878462327292 & -1.8347489733802302 &
   1.5708531406471731 & 2.7E-12 & 6.00000 \\   
  (1,0) $(1,+,-,-)$   & 2.0813680294942181  &   -2.5698934824440314  & 
  -0.49109532951750123 & 
  4.7125297893702216 & 1.2E-13 & 6.00862 \\
  $(1,+,-,+)$ & 2.0813650462038473   &    -2.5665777952951938  & 0.49437520287388081  & 4.7125274954241636  &   4.8E-12  & 6.00698   \\
  (2,0) $(1,+,+,+)$ &  3.6603651202172944 &      -3.2773693881342307 &  0.35681229636585665 &
  10.995511664216785 & 3.6E-13 & 6.00001  \\
   $(1,+,-,-)$ &  2.9251498031208745 &      -3.3408470966236936 & -0.41036694348536390 
   &   7.8540763304982182  & 2.2E-13 & 6.00020 \\
  $(1,+,-,+)$ & 3.6604430591190056 &      -3.1643788287931724 & -0.16773314906446030      
   & 10.995485072271800  & 4.0E-13 &  6.00000 \\
  $(1,+,+,-)$ &  2.9251480066166442   &    -3.3381414081774690 &    0.41309419532725156 &
   7.8540754560494070 & 7.5E-12 & 6.00021 \\
   $(1,-,-,+)$ & -33.386021747563461 &       33.213909690557166  & -2.71073E-015  &
   9.4738871457265912  & 2.0E-11 & 6.00000 \\     
   (3,0) $(1,-,-,-)$ &  3.6605020710174028 &      -4.1787039671068289  &  6.1282794475096E-002 &
   10.995670176360608 & 4.8E-12 & 6.00000  \\
   $(1,+,-,+)$ & 4.3277544548543272 &      -4.0497829282480664 & -0.39291409418382001 &   14.137133009934892 & 4.9E-12 & 6.00000 
    \\
   $(1,+,+,-)$ & 3.6603792900559191  &     -4.0294477233585688 & 0.36927357467952826 &
   10.995634467872524 & 1.2E-13 & 6.00000 \\       
  (4,0) $(1,+,+,+)$ &  4.9471018385520535 &      -4.6471274202861119 &  0.33580758410801892 & 17.278729051746158 & 8.4E-13 & 6.00000 \\
  $(1,+,-,-)$ & 4.3278004629295896  &     -4.6714061401489113 &  -0.33526306450003002 &
  14.137210486904669 & 3.1E-12 & 6.00000 \\
  $(1,+,-,+)$ & 4.9471069002892705  &     -4.6355614258252587 & -0.32513500590766048 &   17.278727722266336 & 3.3E-13 & 6.00000  \\
  $(1,+,+,-)$ & 4.3277976042220825 &      -4.6662669764705162 & 0.34046217998085415  &   14.137209723765521 &  8.4E-14 & 6.00000  
  \\
   (5,0) $(1,-,-,-)$ &  1.5694476869591147 &     -2.5347623301411661 &  -1.9351494E-13 & 15.707071606841449  & 8.7E-12 & 6.02363 \\
   $(1,+,+,-)$ &  4.9471222846504421  &     -5.2644235905578336 & 0.31763985210876378 &
   17.278792357309545 &  6.00000 \\      
   		\bottomrule
   	\end{tabular}\label{Table3}
   \end{table}

  \begin{table}
   	\centering 
   	\tiny
   	\caption{ A list of initial values of the
   	 Hill-lunar-type of \textsc{SDSP} with
   	$\cos^{2} i=1/2$.}
   	\begin{tabular}{c|llllll}
   		\toprule
(k,j) type  & $\xi_{1}$ &  $\dot{\xi}_{2}$  &  $\dot{\xi}_{3}$   &  $T/4$  & accuracy & $\rho$  \\
   		\hline
   	(0,4) $(1,+,+,+)$ & 0.19458458234778178 &       1.5950280529591743 & 1.4247117508651472	&
   	1.4965614757267209 & 2.1E-12 & 6.00000 \\
   	$(1,+,+,-)$ & 0.19479555207814894  &      1.5333467594669177 & -1.4967190735481686 &
   	1.4985745486163704 & 7.0E-13 & 6.00000 \\
	(0,5) $(1,+,-,+)$ & 0.20883475231870061 &       -1.7966587251692738 &  1.5312112883077162 &
   		1.6542335677818685 & 2.8E-12 & 6.00000 \\   		
   	(0,6) $(1,+,-,-)$ &	 0.18573890932212353 & -1.8244214661335518 &  -1.6615950540450048&
   		 1.6385282571508664 & 7.1E-13 & 6.00000 \\
   (0,7) $(1,+,-,-)$ & 0.16816807952651075 &      -1.8597785842093075  & -1.7710466713073072 &
   1.6279262696022292 & 1.6E-12 & 6.00000  \\
   (0,8) $(1,+,+,+)$ & 0.13744008315942863 &        2.0202381771564175 & 1.6317319026603057 &
   1.5246240934921413 & 3.5E-13 & 6.00000 \\
   $(1,-,+,+)$ &  -0.15455014961447458 &       2.0450425050414500 & 1.7146972717335558 &
   1.6251745733609688 & 3.9E-12 & 6.00000 \\
   $(1,+,-,+)$ & 0.15437911047248612 &      -1.9598070470369613 & 1.8060620355549017 &
   1.6223019252187139 & 4.9E-12 & 6.00000 \\
  $(1,-,-,+)$ & -0.13765663140958340 &       -1.8448133435768497 & 1.837470575594241 &
  1.5278938496551244 & 4.7E-13 & 6.00000  \\
  (0,9) $(1,+,-,-)$ & 0.14289622406891708  &     -1.9266767540482233  & -1.9632238681425433 &
  1.6141067018642969 & 2.0E-12 & 6.00000 \\
(0,10) $(1,+,+,+)$ & 0.12159837938035424 & 2.0270650757700204  & 1.9076768240912407  & 1.5343386249413951 & 9.4E-13 & 6.00000 \\
 $(1,+,-,-)$ & 0.13360367656266917 & -2.0954815284812112 & -1.9154997962270426 & 1.6126854371306503
  & 1.6E-13 & 6.00000 \\    		
  $(1,+,+,-)$ & 0.12159167039529933 & 2.0348664875498046  &      -1.8989819872184843 & 1.5342210264366785 
  & 4.8E-12 &  6.00000 \\
 $(1,-,-,+)$ & -0.12150492233283411 & -2.1359757477986805 &  1.7793716039307683 & 1.5327029511595240 & 6.3E-12
 & 6.00000  \\ 
 $(1,-,+,-)$ & -0.13345847029345070 & 1.9810562541769794 &  -2.0279836849711224  & 1.6099183076072057 
 &  6.7E-12 & 6.00000 \\  
 (0,20) $(1,+,+,-)$ & 8.07853845833577E-2 & 2.4552287589143869 & -2.4414577005855942 & 1.5515177229234887 & 4.3E-12 & 6.00000  \\ 
 $(+,-,-,-)$  & 0.08481899200869922 & 
       -2.4982470520694178 &  -2.4452069713398923 &   1.5914167828431405 &  2.0E-11 & 6.00000  \\   
   \hline
 (k,j) type  & $\tilde{\xi}_{1}$ &  $\tilde{\xi}_{3}$ &  $\dot{\tilde{\xi}}_{2}$ &  $T/4$ & accuracy & $\rho$  \\
 \hline
   (0,2) $(2,+,-,+)$ & 0.12038642855020419 &     -0.23158072278374456 & 1.8679973545987234 &
   1.5081253549785989 &  3.4E-13 & 6.00000 \\
   (0,3) $(2,+,-,-)$ & 0.19573418852524427 &      -0.20920773071089563 & -2.0895429237612499 &
   1.7123801252911697 &  1.0E-11 & 6.00000 \\    
  (0,5) $(2,+,-,-)$ & 0.14149505268610094  &    -0.15282899985441442 & -2.3449761662964281 &
  1.6535652804050742 & 3.0E-12 & 6.00000  \\ 
 (0,6) $(2,+,+,+)$ &  0.11078036428234198  &     0.11600297935020890 & 2.3950716405810493  &
 1.5217015549202837 & 3.2E-12 & 6.00000  \\    	
   		\bottomrule
   	\end{tabular}\label{Table4}
   \end{table}
   
%
     
   \section{Discussion}
    \label{sec:6}
    
    This paper considers the numerical continuation of the 
    comet- and Hill-type \textsc{DSPO} of
    the approximate system to the full system. The full system can be the \textsc{CRTBP} or Hill's lunar problem. There are few numerical results about Howison and Meyer's \textsc{SDSP}s.
   The numerical results reveal that \textsc{SDSP}s
   exist even if the small
   symplectic scaled parameters 
   are not necessarily too small. 
   Classical continuation method
   based on the implicit function theorem
   which just has a local convergence
   may fail in continuing the approximate 
   initial values to the exact values for
   the \textsc{SDSP}s. 
   Broyden's method with a line-search
   is taken in use to solve the 
   nonlinear periodicity conditions
   expressed in the rectangular coordinates.
   The algorithm and related fortran programs 
   can be referred to \citet{Fortran90}.    
   Numerically, different approximate
   initial values are continued to 
   different initial values for the \textsc{SDSP}s,
   and we use sixteen cases to continue the 
   known initial values of the Keplerian orbits.
   Some cases even fail because of the strong
   perturbation or the essential instability.  
   Many successful examples in Table \ref{Table1},
   \ref{Table3}, \ref{Table4} and in the Figures show that our numerical scheme behaves well in determining
   the \textsc{SDSP}s.
   
   For a periodic orbit with double symmetry,
   only one fourth period information is
   needed to determine the whole orbit, 
   and the reason is well explained.
   We find the characteristic multipliers 
   can be calculated this way. If the orbit starts
   from $\mathscr{L}_{1}^{(0)}$, then
   $\mathscr{R}_{2}\mathcal{Z}_{T/2}\mathscr{R}_{1}
   \mathcal{Z}_{T/2}$ is diagonal with the 
   characteristic multipliers on the diagonal line.
   If the orbit starts from $\mathscr{L}_{2}^{(0)}$,
   the $\mathscr{R}_{1}\tilde{\mathcal{Z}}_{T/2}
   \mathscr{R}_{2}\tilde{\mathcal{Z}}_{T/2}$
   is diagonal, and the characteristic multipliers
   are also on the diagonal line. This proposition
   can be explained geometrically, but there is a lack
   of proof. This work may be considered in future.

    Some questions are raised at the end of this paper. First, numerical method has its advantages and shortcomings. If the perturbation is sufficiently small,
    the period ratio between the inner orbit and the outer orbit is also very
    small and the integration time will be long. If the perturbation is
    mild strong, the stability of the orbit starting from the beginning
    initial values affects the numerical integration and the continuation results.
    The grid search, Poincar\'{e} section method, the optimization method,
    the Lyapunov indicator method and so on can be integrated in order to
    give a better understanding of the three-dimensional phase structure.
    If the energy is fixed, then there are only two unknown initial parameters,
    and the families of periodic orbits can be tracked. The bifurcation of the 
    families will be interesting for analysis. Second, there may exist
    invariant torus around the linearly stable periodic orbits. The linear stability problem
    may be settled with the help of Floquet theory and some perturbation methods. 
    At last, as there exist singularities, so the regularization may be used
     in order to find more spatial periodic orbits in these problems with a background of real astronomy.

\begin{acknowledgements}
The author would like to thank the two anonymous reviewers 
of this paper for their constructive comments and suggestions.  
The author is supported by the NSFC, Grant No. 11703006.
\end{acknowledgements}


\bibliographystyle{spr-chicago}
\nocite{*}

\end{document}